\documentclass[twosided, 12pt]{article}

\usepackage{microtype}
\usepackage[cmex10]{amsmath}
\usepackage{amssymb, bm, bbm}
\usepackage{amsthm}		   
\usepackage{enumerate}     
\usepackage{graphicx}      
\usepackage{subfig}        
\usepackage{setspace}
\usepackage{color}
\usepackage[compress]{natbib}
\usepackage{multirow}
\usepackage{rotating}	   
\usepackage[affil-it]{authblk}  
\usepackage{hyperref}
\usepackage[all]{hypcap}   
\usepackage{framed, color}

\newtheorem{theorem}{Theorem}
\newtheorem{lemma}{Lemma}

\newcommand{\matlab}{\textsc{Matlab}\textsuperscript{\textregistered}}

\providecommand{\bs}[1]{\boldsymbol{\mathrm{#1}}}
\providecommand{\abs}[1]{\left\lvert#1\right\rvert}

\providecommand{\E}[1]{{\mathbb E}\left[ #1 \right]}

\providecommand{\cov}[1]{\operatorname{Cov}\left( #1 \right)}

\DeclareMathOperator{\erfc}{erfc}

\newcommand{\JPDF}{f_{\hat{R}, \hat{\theta}}(R,\theta)}
\newcommand{\NPJPDF}{f^{(0)}_{\hat{R}, \hat{\theta}}(R,\theta)}
\newcommand{\bfs}{\bs{s}}
\newcommand{\bfr}{\bs{r}}
\newcommand{\rr}{\|\bs{r}\|}
\newcommand{\bfk}{\bs{k}}
\newcommand{\kk}{\|\bs{k}\|}

\newcommand{\stdx}{\sigma}

\newcommand{\hR}{\hat{R}}
\newcommand{\eR}{{R}^{\ast}}
\newcommand{\cR}{\check{R}}
\newcommand{\htheta}{\hat{\theta}}
\newcommand{\ethe}{{\theta}^{\ast}}
\newcommand{\ctheta}{\check{\theta}}
\newcommand{\cxx}{c}
\newcommand{\mx}{m}
\newcommand{\Q}{{\mathbb{Q}}}
\newcommand{\Qind}{\hat{\mathbb{Q}}}
\newcommand{\CQQ}{\bs{C}_{\Qind}}

\newcommand{\hQ}{\hat{Q}}
\newcommand{\eQ}{{Q}^{\ast}}
\newcommand{\cQ}{\check{Q}}
\newcommand{\qd}{q_{\mathrm d}}
\newcommand{\qo}{q_{\mathrm o}}
\newcommand{\eqd}{q_{\mathrm d}^{\ast}}
\newcommand{\eqo}{q_{\mathrm o}^{\ast}}

\newcommand{\bq}{\bs{q}}
\newcommand{\hbq}{\bs{\hat{q}}}
\newcommand{\Rd}{\mathbb{R}}
\newcommand{\ee}{\mathrm {e}}
\newcommand{\ra}{\rightarrow}
\newcommand{\la}{\lambda}

\begin{document}

\title{Non-Parametric Approximations for Anisotropy Estimation in Two-dimensional Differentiable Gaussian Random Fields}

\author{Manolis P. Petrakis%
	\thanks{\href{mailto:petrakis@mred.tuc.gr}{\texttt{petrakis@mred.tuc.gr}}} }
\author{Dionissios T. Hristopulos%
	\thanks{Corresponding author; \href{mailto:dionisi@mred.tuc.gr}{\texttt{dionisi@mred.tuc.gr}}}}
\affil{Geostatistics Laboratory, School of Mineral Resources Engineering, Technical University of Crete, Chania 73100, Greece}

\date{Dated: \today}

\maketitle

\begin{abstract}
    Spatially referenced data often have autocovariance functions with elliptical isolevel contours,
    a property known as geometric anisotropy. The anisotropy parameters include the tilt of the ellipse (orientation angle) with respect to a reference axis and the aspect ratio of the principal correlation lengths. Since these parameters are unknown \emph{a priori}, sample estimates are needed to define suitable spatial models for the interpolation of incomplete data. The distribution of the anisotropy statistics is determined by a non-Gaussian sampling joint probability density. By means of analytical calculations, we derive an explicit expression for the joint probability density function of the anisotropy statistics for Gaussian, stationary and differentiable random fields. Based on this expression, we obtain an approximate joint density which we use to formulate a statistical test for isotropy. The approximate joint density is independent of the autocovariance function and provides conservative probability and confidence regions for  the anisotropy parameters.  We validate the theoretical analysis by means of simulations using synthetic data, and we illustrate the detection of anisotropy changes with a case study involving background radiation exposure data. The approximate joint density provides (i) a stand-alone approximate estimate of the anisotropy statistics distribution (ii) informed initial values for maximum likelihood estimation,
    and (iii) a useful prior for Bayesian anisotropy inference.
    
\end{abstract}


\section{Introduction}

Fast and accurate methods of anisotropy estimation are needed in various fields to better model spatially extended processes and  the properties of heterogeneous materials~\citep{Guilleminot12}. The characterization and measurement of anisotropy in biological tissues, for example, is important for diagnostic and medical
reasons~\citep{Ranganathan11,Richard10}. Significant changes in anisotropy over time may suggest a crucial change in the underlying physical processes. For example, an accidental release of radioactivity may significantly alter the anisotropy of radioactivity patterns over the monitored area. Reliable and computationally fast detection of systematic changes in spatial distributions is crucial, especially for automatic monitoring systems~\citep{INTAMAP}. Another practical question is what constitutes a significant departure from isotropy to necessitate the use of anisotropic autocovariance functions. Non-parametric methods attempt to provide answers to such questions without requiring knowledge of the autocovariance functions (henceforward, \emph{covariance function} for simplicity). Non-parametric isotropy tests are thoroughly reviewed in a recent publication~\citep{Weller15}.

Two types of anisotropy are usually encountered in spatially extended processes. \emph{Physical anisotropy} implies tensor fields that represent directionally dependent material properties such as transport coefficients in heterogeneous media, e.g.~\citep{Adler92}. \emph{Statistical anisotropy} characterizes scalar processes (e.g., scalar permeability, pollutant concentrations), the correlation range of which depends on the spatial direction. Geostatistical analysis employs two types of statistical anisotropy: geometric and zonal~\citep{Zimmerman93,Chiles12}. Herein we focus on geometric anisotropy, which implies SRFs with covariance functions that possess elliptical isolevel contours (see Fig.~\ref{fig:anisotropy_ellipse}). The estimation of anisotropy parameters is a topic of ongoing interest in various engineering fields~\citep{Jiang05,Okada05,Feng08,Olhede08,Bihan01,Xu09,Richard10,Wang2011} and in data assimilation~\citep{Weaver2012diffusion}. In geostatistics, the anisotropy is typically modeled by estimating the empirical variogram in different directions and fitting anisotropic variogram models~\citep{Chiles12}. For second-order stationary SRFs the variogram is equivalent to the covariance function. However, the interpretation of such variogram analysis is not always straightforward~\citep{Weller15}. Anisotropic modeling in the Bayesian framework has also been investigated~\citep{Ecker99,Ecker03}. Recently, there is interest in anisotropic models with locally varying parameters~\citep{Lillah15}. A study focusing on general characterizations of anisotropy beyond the geometrical model appears in~\citep{Allard15}.

The mathematical framework for the study of anisotropy in spatial processes is based on spatial random fields (SRFs), also known as spatial random functions~\citep{Adler81,Christakos1992random,Wackernagel03,Lantu02}. SRFs are used in several scientific and engineering disciplines that study spatially distributed processes (e.g., image processing, theory of transport in heterogeneous media, wave propagation in random media, environmental modeling). SRFs with Gaussian joint probability density function also provide the mathematical framework of \emph{Gaussian processes} in machine learning. Spatially referenced data are typically modeled as SRFs. The analysis of SRFs based on data involves a number of distributional assumptions that need to be validated. A common assumption is that of statistical stationarity which states that the statistical properties are independent of the position. The less strict second-order stationarity assumption is used in practice and requires the expectation of the field to be constant and the covariance function to depend only on the spatial lag. In the case of Gaussian random fields, second-order stationarity is equivalent to strong stationarity. Isotropy is a stricter assumption that requires the covariance function to depend only on the magnitude but not on the direction of the lag. For convenience, isotropic SRF models are often used, even though many real data sets display anisotropic patterns.

In the case of two-dimensional SRFs that admit first-order derivatives in the mean-square sense, a non-parametric and non-iterative method for semi-analytic estimation of anisotropy parameters was proposed and studied in~\citep{dth02,dth08}. This manuscript extends the works above by investigating the joint dependence of the anisotropy parameter estimates. We derive a non-parametric approximation of the \emph{sampling joint probability density function (JPDF)} of anisotropy statistics for differentiable, stationary Gaussian SRFs. We prove this expression using the Covariance Hessian Identity (CHI)~\citep{Swerling62}, the Central Limit Theorem, Jacobi's multivariate  transformation theorem, and perturbation analysis.

The non-parametric approximation yields a sampling JPDF which is more dispersed in parameter space than the exact JPDF. This implies wider probability regions for the anisotropy parameter statistics and confidence regions for the anisotropy population parameters. Hence, if a sample is classified as isotropic at confidence level $p$ based on the approximate JPDF, it is actually isotropic at $p'>p$. The JPDF that we derive can also be used as a prior in Bayesian model inference~\citep{Ecker99,Ecker03,Schmidt03,Zhang12} or as a preliminary step in copula-based spatial analysis~\citep{Kazianka13}.

This manuscript is structured as follows: In Section~\ref{sec:CHI} we present essential definitions and an overview of CHI. In Section~\ref{sec:joint-pdf} we derive a general expression for the joint probability density $\JPDF$ for the anisotropy statistics $(\hat{R},\hat{\theta})$. In addition, we obtain a relation for $p$-level probability regions of the anisotropy parameters. In Section~\ref{sec:non-parametric}, we derive the non-parametric approximation of $\JPDF$ and the corresponding probability region expression. In Section~\ref{sec:isot-test} we formulate a non-parametric test for isotropy. In Section~\ref{sec:simulations}, we validate the theoretical results with numerical simulations and we illustrate the detection of anisotropy changes with a case study involving radiation exposure data. Finally, in Section~\ref{sec:concl} we review the main results obtained in this work, we present our conclusions, and we outline directions for future research. Proofs of theorems and lemmas are given in the Appendices.


\section{Preliminaries}
\label{sec:CHI}
We use boldface symbols for vectors, matrices and tensors; the superscript ``t'' denotes the vector or matrix \emph{transpose}. ${\cal D} \subset \mathbb{R}^{2}$ denotes the spatial domain, $|{\cal D}|$ the enclosed area, $\bfs \in {\cal D}$ the position vector in ${\cal D}$, and $\|\bfs\|$ the Euclidean norm of $\bfs$. $X(\bfs,\omega)$ represents a scalar SRF on the probability space $(\Omega ,\mathcal{F,P})$. The state index $\omega$ determines the field state and is suppressed in the following for the sake of brevity. The events in $\mathcal{F}$ comprise the measured SRF realization(s) or \emph{sample state(s)}. $\E{ \cdot }$ denotes the expectation over the ensemble of states, and $\cov{Z_{1},Z_{2}}= \E{Z_{1}Z_{2}} - \E{Z_{1}}\E{Z_{2}}$  is the \emph{covariance} of the random variables $Z_{1}$ and $Z_{2}$. Realizations of an SRF $X(\bs{s})$ will be denoted by $x(\bs{s})$.

We focus on \emph{wide-sense stationary} Gaussian SRFs (GSRFs) with constant mean $\mx= \E{ X(\bfs) }$ and \textit{covariance function} $\cxx(\bfr)= \E{ X(\bfs) \, X(\bfs +\bfr)} - \mx^{2}$. We assume that the SRF is first-order differentiable in the mean square sense for every $\bfs\in {\cal D}$, so that the partial derivatives ${\partial^2 \cxx (\bfr)} /{\partial r_i^2 }$ in the orthogonal directions $i = 1,2$ exist at $\bfr = (0,0)$. For Gaussian SRFs, mean square differentiability essentially implies that the respective
derivatives of the sample states exist almost surely~\citep{Adler81,Yaglom87}. We assume \emph{short-range} correlations, i.e., with a finite \emph{correlation area} $\int d{\bf r} \, |\cxx({\bf r})|$. Such correlation functions have a finite integral range.

\begin{figure}
	\centering
	\includegraphics[width=0.5\linewidth]{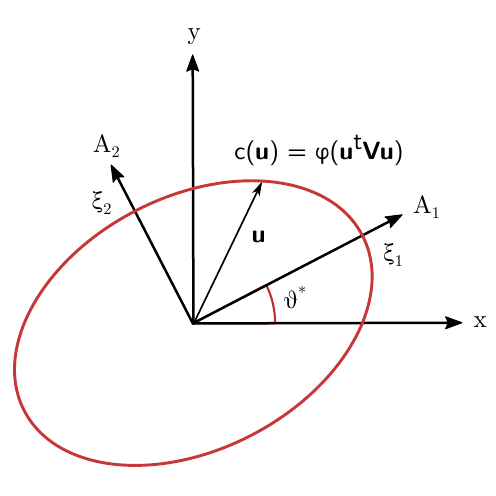}
	\caption{Definition of elliptical anisotropy parameters. The ellipse denotes an iso-level contour of an anisotropic covariance function $c(\cdot)$. The principal correlation lengths are $\xi_1, \xi_2$; $R^\ast=\xi_2/\xi_1$ is the anisotropy ratio, while the $x$-axis and $A_1$ are the sides of the anisotropy angle $\theta^\ast$. An anisotropic covariance function $c(\cdot)$ is obtained from a positive definite function $\varphi(\cdot)$ via a rescaling $\bs{V}$ followed by a rotation by $\theta^\ast$.}
	\label{fig:anisotropy_ellipse}
\end{figure}

The \emph{sample}, ${\bf x}_{k}=(x_{1}, \ldots, x_{N})^{t}$ comprises the values $x_{k} = x(\bs{s}_{k})$ of the realization $x(\bs{s})$, where $\bs{s}_{k}, \, k=1, \dots, N$ are sampling locations. We use the following notation for the anisotropic parameters, illustrated  in terms of the anisotropic ratio: population parameters are marked by a star, i.e., $\eR$. The sampling function of $\eR$ is the random variable $\hR$. Specific numerical values will be denoted by $R$. Sampling functions based on discrete approximations of derivatives are denoted by $\cR$. The population anisotropy parameters are illustrated in Fig.~\ref{fig:anisotropy_ellipse}.

The \textit{Covariance Hessian Matrix} $\bs{H}(\bs{r})$ (CHM) of a stationary, at least first-order differentiable, SRF $X(\bfs)$ is defined as follows
\begin{equation}
    \label{eq:H}
    H_{ij} (\bfr) \doteq -\frac{\partial ^2 \cxx(\bfr)} {\partial r_i \, \partial r_j }, \quad i, j = 1, 2 .
\end{equation}

Let $X_{ij}(\bfs) = \partial_{i} X(\bfs) \, \partial_{j} X(\bfs)$, $i=1,2$ be the \textit{gradient tensor}, where $\partial_{i} X(\bfs)=\partial_{i} X(\bfs)/\partial s_{i}$, $i=1,2$ are the partial derivatives of $X(\bs{s})$. The mean gradient tensor $\bs{\eQ}$, also known as the matrix of spectral moments~\citep{Adler81}, is defined as follows
\begin{equation}
\label{eq:exp-qij}
    \eQ_{ij} \doteq \E{ \partial_{i} X(\bfs) \, \partial_{j} X(\bfs) } = \E{X_{ij}(\bfs)}.
\end{equation}
The matrix $\bs{\eQ}$ is nonnegative definite as the covariance of the random gradient $\nabla X(\bfs)= \left(\partial_1 X(\bs{s}), \partial_2 X(\bs{s})\right)^t$. It satisfies the following theorem~\citep{Swerling62}:

\begin{theorem}[Swerling's CHI]
    \label{theor:chi}
    Let $X(\bfs)$ be a statistically stationary SRF with  covariance function $\cxx (\bfr)$ that admits
    partial derivatives ${\partial^2 \cxx (\bfr)} /{\partial r_i^2 }$ at $\bfr = (0,0)$. Then
    \begin{equation}
        \label{chi}
        \bs{\eQ} = \left. \bs{H} (\bs{r}) \right|_{\bs{r} = \bs{0}}.
    \end{equation}
\end{theorem}

To define the anisotropy parameters, consider a coordinate system aligned with the principal axes of anisotropy, e.g., $A_1$ and $A_2$ (see Fig.~\ref{fig:anisotropy_ellipse}). In the principal system, $\cxx({\bf u})= \phi({\bf u}^{t} {\bf V} {\bf u})$, where ${\bf u}=(u_1, u_2)$ is the lag, $\bs{V}$ is a diagonal $2 \times 2$ matrix, and $\varphi(\cdot)$ is a positive definite function.
The principal correlation lengths of $X({\bfs})$ are given by $\xi^{-2}_{i} = -a \stdx^{-2} \, \partial^{2} \cxx({\bf u})/\partial {u}_i^{2} |_{\bf u={\bf 0}}$, for $i=1,2$ where $a$ is a positive constant~\citep{dth08}. We define the anisotropy ratio as $\eR=\xi_2/\xi_1 $, and the orientation (rotation) angle $\ethe$ as the angle between the horizontal axis of the reference system and $A_1$. The anisotropy parameters $(\eR,\ethe)$ satisfy the following theorem:

\begin{theorem}
    \label{theor:aniso}
    Let $X({\bfs})$ be an SRF satisfying the conditions of Theorem~\ref{theor:chi}, and $q_{\rm d}$, $q_{\rm o}$
    represent the following ratios of gradient tensor elements $\eQ_{ij}$:
    \begin{subequations}
    \label{eq:anisotropy}
    \begin{equation}
        \label{eq:qdiag}
        \eqd \doteq \frac{\eQ_{22} }{\eQ_{11} }= {\frac{1+ {\eR}^2 \tan^2\ethe }{{\eR}^2 +\tan^2\ethe}},
    \end{equation}
    \begin{equation}
        \label{eq:qoff}
        \eqo \doteq \frac{\eQ_{12} }{\eQ_{11} }=\frac{\tan \ethe ({\eR}^2 - 1)}{{\eR}^2 + \tan^2\ethe}.
    \end{equation}
    \end{subequations}

    Then, the anisotropic ratio, $\eR$ and the orientation angle, $\ethe$ are given by
    \begin{subequations}
    \label{eq:R-theta}
    \begin{align}
        \ethe & =    \frac{1}{2}\tan^{-1}\bigg(\frac{2 \eqo}{1- \eqd}\bigg),
        \label{eq:theta-sol}    \\
        \eR &   =   \left[1+\frac{1 - \eqd}{ \eqd -(1+ \eqd)\cos^2\ethe}\right]^{-1/2}.
        \label{eq:R-sol}
    \end{align}
    \end{subequations}
\end{theorem}

The proof is based on Theorem~\ref{theor:chi}~\citep{dth08}. Therein the notation $R=R_{2(1)}=\xi_1/\xi_2$ was used, whereas above we defined $\eR = \xi_2/\xi_1$. The equations~\eqref{eq:anisotropy} and~\eqref{eq:R-theta} follow from~\citep{dth08} by means of the transformation $R \mapsto 1/\eR$.

Equations~\eqref{eq:anisotropy} are invariant under the pair of transformations $\tan\ethe \mapsto - (\tan \ethe)^{-1}$, that is, $\ethe \mapsto \ethe \pm \pi/2,$ and $\eR \mapsto 1/\eR$. By restricting the parameter space to $\eR \in [0,\infty)$ and $\ethe \in [-\pi/4,\pi/4)$, or equivalently to $\eR \in [1,\infty)$ and $\ethe \in [-\pi/2,\pi/2)$, ensures that the mapping $(\eqd,\eqo) \mapsto (\eR,\ethe)$ is one-to-one, except for the point $(1,0)$ which maps to $(1,\ethe)$  $\ethe$ being any angle $\in [-\pi/2,\pi/2)$. Theorem~\ref{theor:aniso} permits estimating the anisotropy parameters without knowledge of the covariance function, if $\bs{\eQ}$ can be estimated from the data~\citep{dth02,dth08}.


\section{Sampling Joint PDF of Anisotropy Statistics}
\label{sec:joint-pdf}

Every realization $x(\bfs)$ yields a different estimate of $\bs{\eQ}$, leading to a probability distribution for the statistics $\hR$ and $\htheta$. Below we derive the joint PDF $\JPDF$ based on Jacobi's theorems for the transformation of a multivariate probability distribution under transformation of the respective variables~\citep{Papoulis} and the Central Limit Theorem (CLT).

We estimate ${\eQ}_{ij}$ using the spatially averaged gradient tensor $\hQ_{ij}$, where  $i,j=1,2$,
\begin{equation}
\label{eq:Q}
    \hQ_{ij} := \frac{1}{N}\sum_{k=1}^{N} {X}_{ij}(\bs{s}_k) =
    \frac{1}{N}\sum_{k=1}^{N} {\partial}_{i}{X}(\bs{s}_k) \,    {\partial}_{j}{X}(\bs{s}_k).
\end{equation}
The estimation of the field's partial derivatives from the data is discussed in Section~\ref{sec:simulations}.
Replacing the expectation with the spatial average requires the ergodic hypothesis. A necessary condition for ergodicity is that $|\mathcal{D}| \to \infty$ in such a way that both ratios of domain length over the correlation length in the respective direction tend to infinity. In practice, this means that for an accurate estimate of ${\eQ}_{ij}$ the domain length along each principal direction should be considerably larger than the
respective correlation length. In the following, we assume that the \emph{asymptotic regime} is defined by
$|\mathcal{D}| \to \infty$ in the sense defined above for ergodicity and $N \to \infty$ (for application of the CLT).


\subsection{Joint PDF of Gradient Tensor Components}

We define the following random vector
\begin{equation}
\label{eq:Qdef}
	{\Qind} = (\hat{Q}_{11},\hat{Q}_{22},\hat{Q}_{12})^t  =
	\left(  \frac{1}{N}\sum_{k=1}^{N} {X}_{11}(\bs{s}_k), \,  \frac{1}{N}\sum_{k=1}^{N} {X}_{22}(\bs{s}_k) , \,
	\frac{1}{N}\sum_{k=1}^{N} {X}_{12}(\bs{s}_k)   \right)^t,
\end{equation}
that comprises the independent components of the \emph{fully symmetric gradient tensor sampling function}
$(\hat{Q}_{12}= \hat{Q}_{21})$. As we show below, $\Qind$ tends to follow the joint Gaussian distribution in the asymptotic limit due to the Central limit theorem.

According to~\eqref{eq:Q}, $\hQ_{ij} = \frac{1}{N} \sum_{k=1}^{N} X_{ij}(\bs{s}_k)$ and based on the definition~\eqref{eq:exp-qij} it follows that $\E{\hQ_{ij}} = \eQ_{ij}$, i.e., $\Qind$ is an unbiased estimator of $\bs{\eQ}$. By definition, the covariance matrix $\CQQ$ is symmetric, namely $C_{ij;kl}=C_{kl,ij}$; hence, it involves six independent entries.

\begin{lemma}[Covariance matrix $\CQQ$]
\label{theor:CQQ}
	For a statistically stationary GSRF, the six independent entries of $\CQQ$ are given by the following series
	\begin{align}
	\label{eq:Cijkl}
	    C_{ij;kl} & = \frac{1}{N^2} \sum_{\bs{r}_{nm}} C_{ij;kl}(\bs{r_{nm}}) =
	    \frac{1}{N^2} \sum_{\bs{r}_{nm}} \left[ H_{ik}(\bs{r}_{nm}) H_{jl}(\bs{r}_{nm}) +
	    H_{il}(\bs{r}_{nm}) H_{jk}(\bs{r}_{nm}) \right]\! \nonumber \\
	    & = \frac{1}{N}  \left[ \eQ_{ik} \, \eQ_{jl} + \eQ_{il} \, \eQ_{jk} \right]  \nonumber \\
	    & + \frac{1}{N^2} \sum_{\bs{r}_{nm} \neq \bs{0}} \left[H_{ik}(\bs{r}_{nm}) H_{jl}(\bs{r}_{nm}) +
	        H_{il}(\bs{r}_{nm}) H_{jk}(\bs{r}_{nm}) \right],
\end{align}
for
\[(i,j,k,l) \in \{(1,1,1,1), (1,1,2,2), (1,2,1,2), (1,1,1,2), (2,2,2,2), (1,2,2,2)\}\]
where $\bs{r}_{nm}=\bs{s}_n - \bs{s}_m$  is the lag vector between two locations $\bs{s}_{n}$ and $\bs{s}_{m}$ for $n, m=1,\dots, N$.
\end{lemma}
\begin{proof}
    The proof is given in Appendix~A. This is the only step in which we employ the Gaussian assumption for the joint PDF in order to accomplish the decomposition of higher than second-order moments based on the Wick-Isserlis theorem. However, the Gaussian assumption could be relaxed using a variational Gaussian approximation.
\end{proof}

The term $N^{-1} \left( \eQ_{ik} \, \eQ_{jl} + \eQ_{il} \, \eQ_{jk} \right)$ in~\eqref{eq:Cijkl}
is obtained from the summands with $\bs{r}_{nm} = \bs{0}$ and leads to the non-parametric approximation of  $\JPDF$ as shown below. The sums over $\bs{r}_{nm}\neq \bs{0}$ include parametric corrections that depend on the covariance function. In the  approximate, non-parametric expression we omit the parametric terms which are smaller. These terms have an $1/N^2$ prefactor, but they also involve  $N^2$ summands. However, the products of the covariance Hessians $H_{ik}(\bs{r}_{nm})\, H_{jl}(\bs{r}_{nm})$ that appear in the summands decay very fast with $\|\bs{r}_{nm}\|$. This is due to the fact that, according to~\eqref{eq:H}, the covariance Hessian decays in space proportionally to the second derivative of the covariance function; assuming ergodic conditions, this decay is fast. Hence, at most $\mathcal{O}(N)$ of these terms, for which $\|\bs{r}_{nm}\| < \min(\xi_{1}, \xi_{2})$, make a significant contribution. Thus, the parametric correction is at most $\mathcal{O}(1/N)$. On the other hand, since it involves $H_{ik}(\bs{r}_{nm})$  at finite lag distances, the corrections are smaller (in absolute value) than the non-parametric component. In the isotropic case, $H_{12}(\bs{r}_{nm})=0$ for every $\bs{r}_{nm}$.

\begin{lemma}[Joint PDF of $\Qind$]
	\label{theor:fQ}
    Assume $X(\bs{s})$ is a statistically stationary SRF  with  short-ranged covariance $\cxx({\bf r})$ whose spectral density satisfies $\tilde{C}({\bfk}) \sim \mathcal{O}(\kk^{-3-\epsilon})$  for $\epsilon>0$ as $\kk \ra \infty$. Then, the joint PDF of the vector $\Qind$ which is defined by~\eqref{eq:Qdef} tends asymptotically to the following \textit{trivariate Gaussian}
    \begin{equation}
        \label{eq:f-v}
        f_{\Qind}(\Qind;\bs{\eQ},{\CQQ}) =  \frac{\ee^{-\frac{1}{2} (\Q - \bs{\eQ})^{t} {\CQQ}^{-1}(\Q - \bs{\eQ})}}{(2\pi)^{3/2} \, \det({\CQQ})^{1/2}},
    \end{equation}
    where $\E{\Qind}= \bs{\eQ}$ and the covariance matrix $\CQQ$ is defined by~\eqref{eq:Cijkl}.
\end{lemma}

\begin{proof}
    The proof is given in the Appendix~B.
    The condition $\tilde{C}({\bfk}) \sim \linebreak[4] \mathcal{O}(\kk^{-3-\epsilon})$, $\epsilon>0$ implies that for every  $ \kk \ra \infty$, there are $ \epsilon >0$ and $ C_{\infty}>0$, such that $\tilde{C}({\bfk}) \le C_{\infty}/\kk^{3+\epsilon}.$ This is satisfied by most finite-range, twice differentiable covariance functions, including the Gaussian, rational quadratic, Bessel-J, and Mat\'{e}rn with $\nu> 1$ covariance models~\citep{Lantu02}.
\end{proof}


\subsection{PDF of Gradient Tensor Ratios}
Based on the joint PDF of $\Qind$, we derive the JPDF of the gradient tensor ratios $f_{\hbq}(\bq;\bs{\eQ},{\CQQ})$, where $\bq =(\qd, \qo)^{t}$.

\begin{lemma}[PDF of gradient tensor ratios]
	\label{theor:fq}
	For an SRF $X(\bs{s}) $ that satisfies the conditions of Lemma~\ref{theor:fQ}, the joint density  $f_{\hbq}(\bq;\bs{\eQ},{\CQQ})$ tends asymptotically to  the following non-Gaussian density
	\begin{equation}
	\label{eq:f-qd-qo}
	    f_{\hbq}(\bq; \bs{\eQ},{\CQQ}) = \frac{\la_{2} \, \ee^{-\frac{\la_{1}}{2}}}{8 z_{1}^{5} } \,
	    \bigg[ \sqrt{2\pi} \, (z_{2}^2 + 4 z_{1}^2) \, \exp\left(\frac{z_{2}^2}{8 z_{1}^{2}}\right) \,
	    \erfc\left(\frac{z_{2}}{2 \sqrt{2}z_{1}}\right) - 4 z_{1}z_{2} \bigg],
	\end{equation}
	
	\noindent where $\erfc(\cdot)$ is the \textit{complementary error function}, and $z_{1}, z_{2}$,
	$\la_{1}, \la_{2}$ in the above expression are given by the following expressions, where ${\bq'}^{t} = (1, \qd, \qo)$
	\begin{subequations}
		\label{eq:triv-pdf-coefs}
		\begin{align}
			z_{1}^{2}(\bq;\CQQ) &=  {\bq'}^{t} \, {\CQQ^{-1}} \, {\bq'},\label{eq:A}\\
			z_{2}(\bq;\bs{\eQ}, {\CQQ}) &= -2  \, \bs{\eQ}^{t} \,  {\CQQ^{-1}} \, \bq', \label{eq:B}\\
			\la_{1}(\bs{\eQ}, {\CQQ}) & = \bs{\eQ}^{t} \, {\CQQ^{-1}} \, \bs{\eQ},\label{eq:C}\\
			\la_{2}({\CQQ}) & = (2 \pi)^{-3/2} \, [\det (\CQQ ) ]^{-1/2}.\label{eq:K}
		\end{align}
	\end{subequations}
\end{lemma}

\begin{proof}
    The proof is based on the transformation of the JPDF under the change of variables $\Qind \mapsto \bs{\hat{q}}$ and is given in Appendix~C.
\end{proof}

We simplify~\eqref{eq:f-qd-qo} by explicitly showing the dependence of $f_{\hbq}(\bq;\bs{\eQ},{\CQQ})$ on $N$.
First, note that as shown by~\eqref{eq:Cijkl} and the associated dimensional analysis, $\CQQ \propto \mathcal{O}(1/N)$. Since $z_{1}^{2}>0$ for all correlated SRFs, we can define $ y\, \sqrt{N} = z_{2}/ (2\sqrt{2} z_{1})$ and $ 2\, N\,\tilde{\la}_{1} = \la_{1} $. In light of $z_{1}>0$ according to~\eqref{eq:A} and $z_{2}<0$ according to~\eqref{eq:B}, it follows that $y<0$. The JPDF is expressed as follows in terms of $y$
\begin{equation}
\label{eq:triv-pdf-scaled}
    f_{\hbq}(\bq; \bs{\eQ},{\CQQ}) =
    \frac{\la_{2} \ee^{-N\,\tilde{\la}_{1}}}{\sqrt{2} z_{1}^{3} } \,
       \Big[\sqrt{\pi} \,(2 y^2 \, N + 1) \, \exp(y^2 \, N) \,
      \erfc\left( y\, \sqrt{N}\right) - 2 \,y\, \sqrt{N}\Big].
\end{equation}

For $y<0$ and $N \to \infty$  we define  $x = y\, \sqrt{N}$, we use the identity $\erfc(x)=2-\erfc(-x)$ and the asymptotic expansion of the complementary error function~\citep[Eq. 7.1.23 and 7.1.24]{Abramowitz1965} to show that
\begin{equation*}
	\erfc(x) = 2 + \ee^{-x^2}\left[\pi^{-1/2} x^{-1}+\mathcal{O}(x^{-2})\right].
\end{equation*}
Hence, to  leading-order in $N$, the JPDF~\eqref{eq:triv-pdf-scaled} is approximated as follows
\begin{equation}
\label{eq:triv-pdf-scaled-approx}
    f_{\hbq}(\bq; \bs{\eQ},{\CQQ}) \approx
       \frac{\sqrt{2\pi} \,\la_{2}  }{ z_{1}^{3} } \,  (2 y^2 \, N + 1) \,
       \exp\left[ (y^2 - \tilde{\la}_{1})\, N\right].
\end{equation}


\subsection{Joint PDF of Anisotropy Statistics}

\begin{theorem}[Joint PDF of anisotropy statistics]
    \label{theor:jpdf}
    For an SRF $X(\bs{s})$ that satisfies the conditions of Lemma~\ref{theor:fQ}, the JPDF of the statistics $\hR$ and $\htheta$ is given asymptotically by
    \begin{equation}
        \label{eq:f-anispar}
        f_{\hR,\htheta}(R,\theta;\bs{\eQ},{\CQQ}) \approx
        \frac{2 R \, | R^2-1 | \,f_{\hbq}(\bq; \bs{\eQ},{\CQQ})}
        {\left(R^2 \cos^2\theta + \sin^2 \theta\right)^3},
    \end{equation}
    where $f_{\hbq}(\bq; \bs{\eQ},{\CQQ})$ is given by~\eqref{eq:triv-pdf-scaled}.
\end{theorem}
\begin{proof}
    The proof is given in Appendix~D. It is based on the transformation of the multivariate probability density
    function under the change of variables  $\bs{q} \mapsto (R, \theta)^t$.
\end{proof}

The function $\JPDF$ is clearly non-Gaussian and depends on $\bs{\eQ}$ and ${\CQQ}$, whereas $\bq$ is expressed in terms of $(R,\theta)$ using~\eqref{eq:qdiag} and~\eqref{eq:qoff}. If the rotation angle is measured in degrees instead of radians, $\JPDF$ should be multiplied by $\pi/180$.


\subsection{Probability Regions for Anisotropy Parameters}
\label{sec:confi-reg}

\begin{figure}
	\centering
	\includegraphics[width=0.85\linewidth]{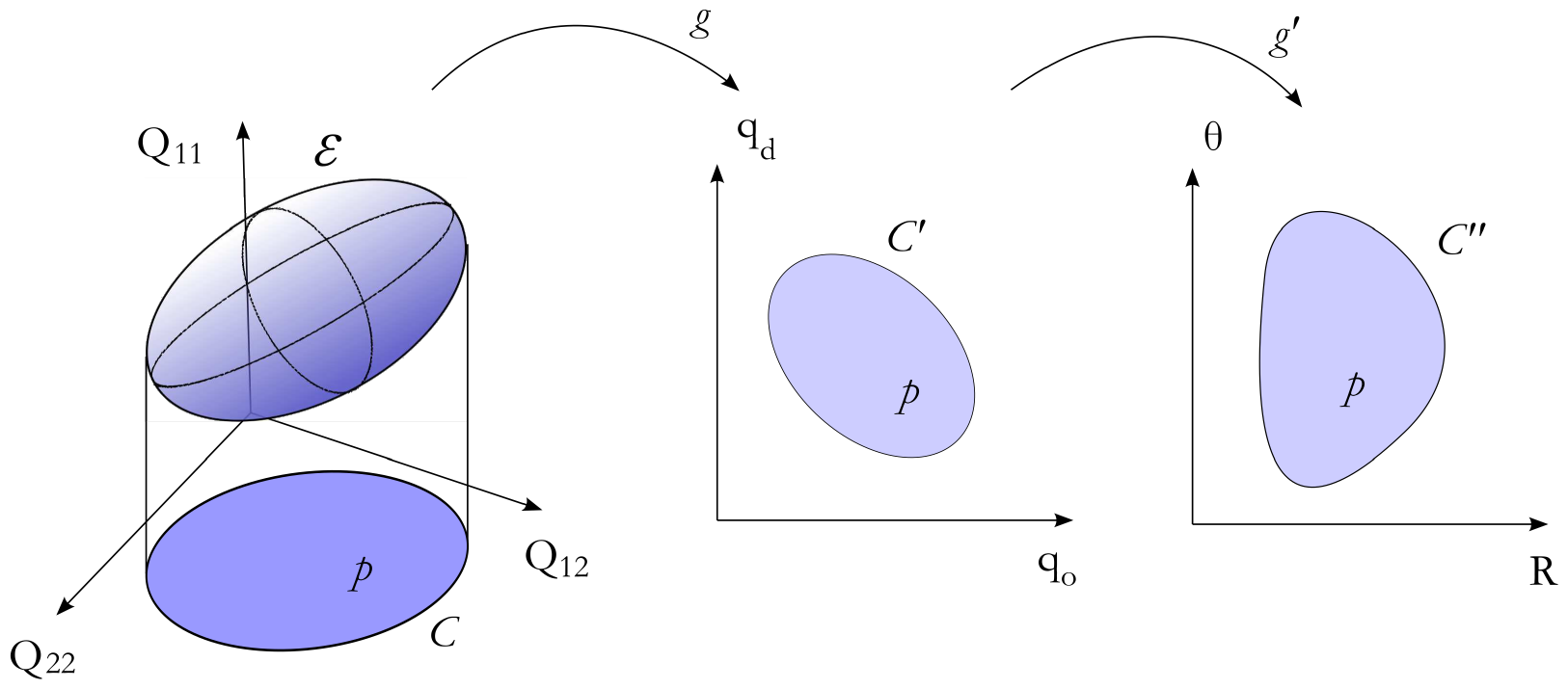}
    \caption{Schematic illustrating the transformation of probability regions at level $p$ due to variable transformations.}
	\label{fig:change_of_variables}
\end{figure}

The probability region at a probability level $p \in [0, 1]$ is the ``volume'' of  space which contains a proportion $p$ of the anisotropy statistics, given the true values $(\eR, \ethe)$.
The probability region of the anisotropy parameters is defined by the following equivalent equations
\begin{equation*}
    p =
    \begin{cases}
        \int_{\mathcal{E}} {\mathrm d}{\Q} \, f_{\Qind}({Q}_{11},{Q}_{22},{Q}_{12};\bs{\eQ},{\CQQ}),\\
        \int_{\mathcal{C}'} {\mathrm d} \qd \, {\mathrm d}\qo \, f_{\hbq}(\qd,\qo;\bs{\eQ},{\CQQ}),\\
         \int_{\mathcal{C}''} {\mathrm d}R \, {\mathrm d}\theta \, f_{\hR,\htheta}(R,\theta;\bs{\eQ},{\CQQ}),
    \end{cases}
\end{equation*}
where $\mathcal{E} \subset \Rd^{3}, \, {\mathcal{C}'} \subset \Rd^{2}$, and $\mathcal{C}'' \subset [0, \infty) \times [-\pi/4,\pi/4)$ represent the probability regions in the respective spaces. These equations represent the evolution of the probability region under the variable transformations $\Q \mapsto \bq \mapsto (R,\theta)$ as shown schematically in Fig.~\ref{fig:change_of_variables}.

\begin{lemma}[Parametric equation of probability regions]
    \label{theor:CR}
    For a SRF $X(\bs{s}) $ that satisfies the conditions of Lemma~\ref{theor:fQ}, the probability region of the anisotropy statistics corresponding to level $p \in [0,1]$ in $(R,\theta)$-space is given by the parametric equation
    \begin{equation}
        \label{eq:confi-reg}
        \left[ y^2(\bq;\bs{\eQ},\CQQ) - \tilde{\la}_{1}(\bq;\bs{\eQ},\CQQ) \right] \,N = \ln(1-p),
    \end{equation}

    \noindent where  $y(\cdot)$, $\tilde{\la}_{1}(\cdot)$ are defined in~\eqref{eq:triv-pdf-scaled} and
    $\bq \mapsto (R,\theta)$ by means of~\eqref{eq:qdiag} and~\eqref{eq:qoff}.
\end{lemma}

\begin{proof}
    The JPDF $f_{\Qind}$ is given by the trivariate Gaussian~\eqref{eq:f-v}. Hence, the probability region of $\Qind$ is an ellipsoid whose surface satisfies the equation
    \begin{equation}
        (\Q - \bs{\eQ})^t \, \CQQ^{-1} \, (\Q - \bs{\eQ}) = \ell_p,
    \end{equation}
    where $\ell_{p}=F^{-1}(\chi^2=p,\nu = 2)$ is the \emph{inverse of the chi-square cumulative distribution
    function} with $\nu=2$ degrees of freedom~\citep{Siotani}. Under the transformation $\Q \mapsto \bq$,
    the ellipsoid is projected onto an ellipse which is deformed by the transformation $\bq \mapsto (R,\theta)$
    into an asymmetric convex curve (see Fig.~\ref{fig:change_of_variables}). Based on \eqref{eq:quadr}, the equation of the corresponding ellipsoid in $(u, \qd, \qo)$-space is given by
    \begin{equation*}
        \label{eq:p-level-qspace}
        z_{1}^{2}(\bq;\CQQ) \, u^2 + z_{2}(\bq;\bs{\eQ}, {\CQQ}) \, u
            + \la_{1}(\bs{\eQ}, {\CQQ}) - \ell_p  = 0,
    \end{equation*}

    \noindent where the coefficients $z_{1}(\cdot), z_{2}(\cdot), \la_{1}(\cdot)$ are given by~\eqref{eq:triv-pdf-coefs}. The above quadratic equation has a unique real solution $u={Q}_{11}$ for any $\bq$ if the discriminant vanishes, i.e.,
    \begin{equation}
        \label{eq:p-level-qdqospace}
        z_{2}^{2}(\bq;\bs{\eQ}, {\CQQ}) -4  z_{1}^{2}(\bq; \CQQ) \, \left[ \la_{1}(\bs{\eQ}, {\CQQ}) - \ell_p\right]=0.
    \end{equation}

    The equation above defines the probability region at level $p$.  We can verify using~\eqref{eq:triv-pdf-coefs} that~\eqref{eq:p-level-qdqospace} represents an ellipse in the space of $\bq$, i.e., it is equivalent to
    \begin{align*}
        \label{eq:discriminant}
        &{\bq'}^{t} \, \bs{M}  \, {\bq'} = 0, \, \text{where} \\
        &\bs{M} = (\CQQ^{-1} \, \bs{\eQ})  \,  (\CQQ^{-1} \, \bs{\eQ})^{t}
          - ( \bs{\eQ}^{t}  \, {\CQQ^{-1}} \, \bs{\eQ} -\ell_{p} )\, {\CQQ^{-1}}.
    \end{align*}

	By inserting in~\eqref{eq:p-level-qdqospace} the functions $y$ and $\tilde{\la}_{1}$ used in~\eqref{eq:triv-pdf-scaled} we obtain the parametric equation $2(\tilde{\la}_{1} - y^2  )\, N  = \ell_p$, where by definition $F(\ell_p, \nu = 2) = p$. Since $F(x,\nu=2)=1-\exp(-x/2)$ \citep[Eq. 26.4.1]{Abramowitz1965}, it follows that $\ell_p = -2\ln(1-p)$, finally leading to~\eqref{eq:confi-reg}.
\end{proof}


\section{Non-parametric JPDF and Probability Region}
\label{sec:non-parametric}

\begin{figure}
    \centering
    \subfloat[Isotropic: $\xi=1$, $\stdx^2=1$.]{\label{subfig:Cijkl_r_iso}
    \includegraphics[width=0.9\linewidth]{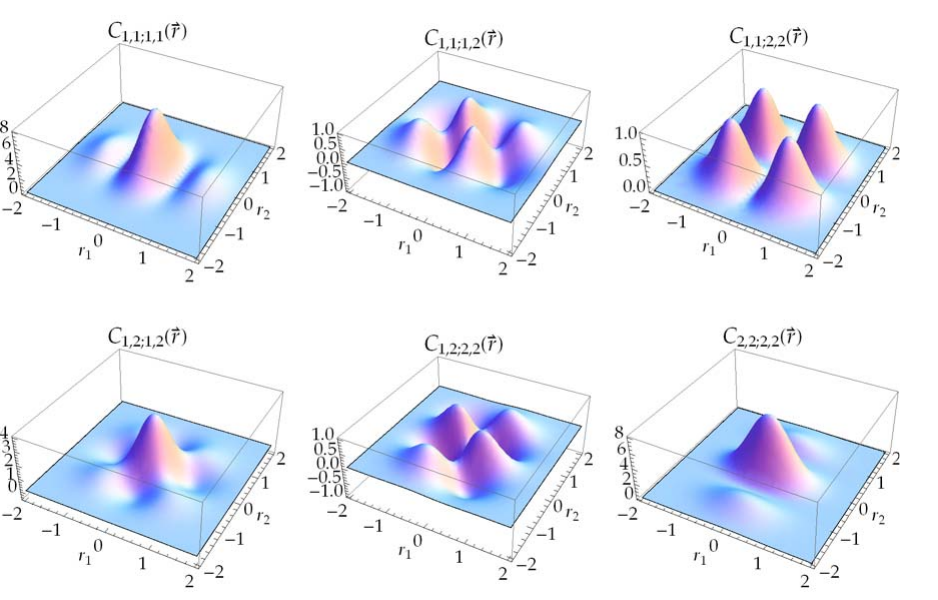}}\\
    \subfloat[Anisotropic: $\xi_1 = 1$, $\xi_2 = 2$, $\theta = 30^\circ$, $\stdx^2=1$.]{\label{subfig:Cijkl_r_aniso}
    \includegraphics[width=0.9\linewidth]{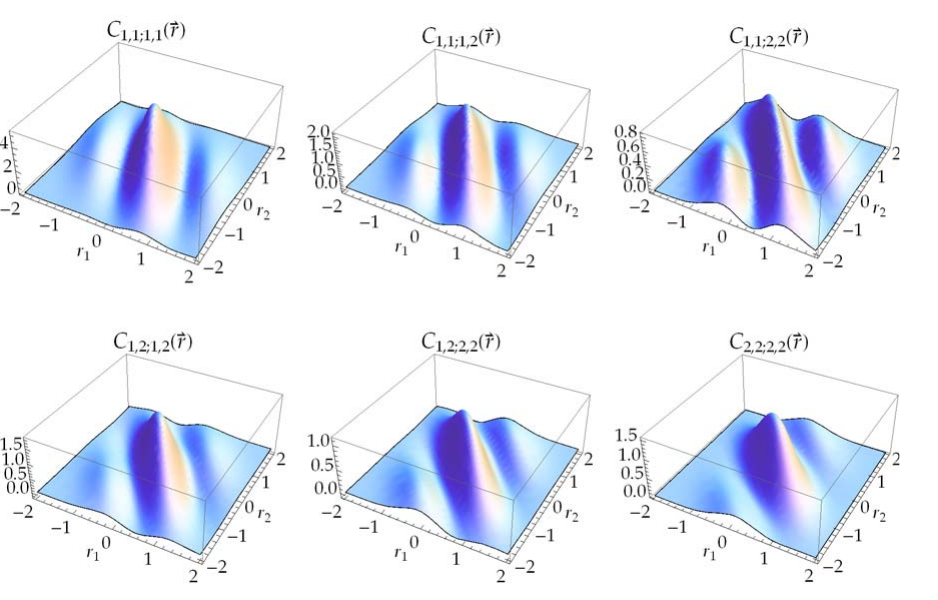}}
    \caption[Plots of $C_{ij;kl}(\bs{r})$ for isotropic and anisotropic Gaussian correlation function]
    {Plots of $C_{ij;kl}(\bs{r})$ for  (a) isotropic $(\xi=1)$ and (b) anisotropic ($\xi_1 = 1$, $\xi_2 = 2$ and $\theta = 30^{\circ}$) Gaussian correlation functions. $C_{ij;kl}(\bs{r})$ essentially vanishes  outside a square of side $a \approx 2\xi$ (isotropic case) and $a \approx 3 \max(\xi_1, \xi_2)$ (anisotropic case).}
    \label{fig:Cijkl_r}
\end{figure}

The expressions for  $f_{\hR,\htheta}(R,\theta;\bs{\eQ},\CQQ)$ and the probability regions of $(\hR,\htheta)$ above depend on the matrix $\CQQ$, given by~\eqref{eq:Cijkl}. $\CQQ$ involves the series~\eqref{eq:Cijkl} that does not, in general,  admit a closed form. If $H_{ij}(\bs{r})$ decays fast for increasing $\|\bs{r}\| $ we can use the explicit approximation  $\CQQ \approx\CQQ^{(0)}$, where
\begin{equation}
    \label{eq:Cvv-nonpar}
    \CQQ^{(0)}= \frac{2}{N} \begin{bmatrix}
    {\eQ_{11}}^2 & {\eQ_{12}}^2 & \eQ_{11}\,\eQ_{12} \\
    {\eQ_{12}}^2 & {\eQ_{22}}^2 & \eQ_{12}\eQ_{22} \\
    \eQ_{11}\eQ_{12} & \eQ_{12}\eQ_{22} & \tfrac{1}{2}({\eQ_{12}}^2 + \eQ_{11} \eQ_{22})
    \end{bmatrix}.
\end{equation}
Figure~\ref{fig:Cijkl_r} illustrates this fast decay of $C_{ij;kl}(\bs{r})$ for isotropic (Fig.~\ref{subfig:Cijkl_r_iso})  and  anisotropic Gaussian covariance (Fig.~\ref{subfig:Cijkl_r_aniso}) functions.
We expect that  $\CQQ^{(0)}$ will lead to a joint PDF with higher uncertainty, and hence more spread out than the true PDF, because it does not incorporate spatial correlations. We validated this intuitive argument by means of numerical simulations (see Section~\ref{ssec:simulations-B}).

\begin{theorem}[Non-parametric JPDF]
	\label{theor:JPDF-NP}
	For an SRF $X(\bs{s})$  that satisfies the conditions of Lemma~\ref{theor:fQ}, the non-parametric JPDF approximation $f^{(0)}_{\hR,\htheta}(R,\theta; \allowbreak \eR,\ethe,N)$ of $(\hR,\htheta)$ is given in the asymptotic regime by
	\begin{subequations}
		\begin{align}
			\label{eq:f-anis-np-jac}
			f^{(0)}_{\hR,\htheta}(R,\theta;\eR,\ethe,N) &  = \abs{\det(\bs{J}_{\theta, R})} f_{\hat{\bs{q}}}^{(0)}(R,\theta;\eR,\ethe,N),
			\intertext{where}
			\label{eq:f-anis-np}
			f_{\hat{\bs{q}}}^{(0)}(R,\theta;\eR,\ethe,N) & \approx \frac{\sqrt{2\pi} \, \la_{2;0}}{z_{1;0}^{3}} \, \left(2 y_{0}^2 \, N + 1\right) \, \ee^{N \left(y_{0}^2 - 1/2 \right)}.
		\end{align}
		\label{eq:f-anis_np_both}
	\end{subequations}
	
	\noindent  The coefficients $z_{1;0}, y_{0}, {\la}_{2;0}$ are given by the following expressions, where $\delta \theta = \theta-\ethe$,
	\begin{subequations}
	\label{eqs:coeff}
	\begin{align}
	    \label{eq:y0}
	    y_{0} &= \frac{1}{ \sqrt{2} \, z_{1;0}}\,\left[(R^2-1)({\eR}^2-1)\cos(2\delta \theta)-(R^2+1)({\eR}^2+1)\right],\\
	    \label{eq:tK}
	    \la_{2;0} &=  \frac{\sqrt{2}\,(\eR)^3}{ \pi^{3/2}} \, \left[(R^2+1)-(R^2-1)\cos(2\theta) \right]^3, \\
	    \label{eq:tA}
	    {z}_{1;0}^{2} &= (R^2-1)^2 ({\eR}^2-1)^2 \cos(4\delta \theta)-4(R^4-1)({\eR}^4-1)\cos(2\delta \theta) \\
	    &   +   (R^4+1)(3 {\eR}^4 + 2 {\eR}^2+3) + 2 R^2 ({\eR}^2-1)^2. \nonumber
	\end{align}
	\end{subequations}
\end{theorem}

\begin{proof}
	In~\eqref{eq:triv-pdf-coefs} we  replace $\CQQ$ with $\CQQ^{(0)}$, defined by~\eqref{eq:Cvv-nonpar}.
	Thus,  $z_{1}, z_{2}, \la_{1}, \la_{2}$ are replaced, respectively, by $z_{1;0}$, $z_{2;0}$, $\la_{1;0}$, $\la_{2;0}$; then, $ y_{0} =z_{2;0} /( 2 \sqrt{N} \, z_{1;0})$. Performing the calculations with $\CQQ^{(0)}$ we obtain~\eqref{eq:y0}-\eqref{eq:tA}. The asymptotic result~\eqref{eq:triv-pdf-scaled-approx} of Lemma~\ref{theor:fq} is used in~\eqref{eq:f-anis-np-jac} to obtain the non-parametric approximation~\eqref{eq:f-anis-np}. Note that in the non-parametric approximation, the coefficient $\tilde{\la}_{1;0}$ in the exponent on the right hand side of~\eqref{eq:f-anis-np} is reduced to $1/2$.
\end{proof}

Numerical comparisons show that the absolute relative error between the non-parametric JPDF $f^{(0)}_{\hR,\htheta}(R,\theta;\eR,\ethe,N)$ calculated with (i) the exact $f_{\hbq}^{(0)}(R,\theta;\eR,\ethe,N)$, obtained from~\eqref{eq:triv-pdf-scaled} by inserting the approximate covariance matrix $\CQQ^{(0)}$, and (ii) the asymptotic limit given by~\eqref{eq:triv-pdf-scaled-approx}, is less than $\approx 10^{-9}$ for $N=50$ and $\approx 10^{-6}$ for $N=30$.

\begin{figure}
    \centering
    \subfloat[$\eR=1.2$, $\ethe=20^\circ$, $N=100$.]%
    {\label{subfig:JPDF_bimodal}\includegraphics[width=0.5\linewidth, trim = 0pt 0pt 0pt 100pt, clip]{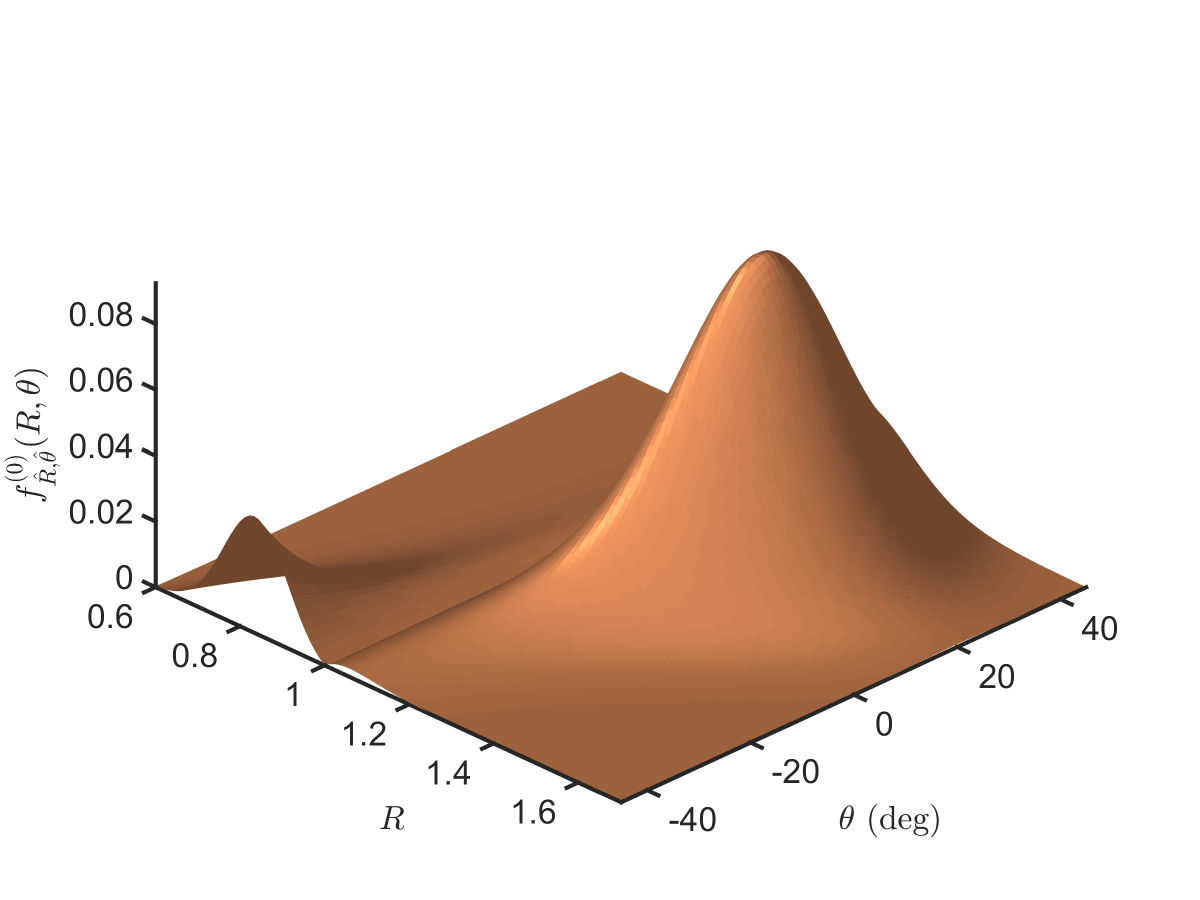}}
    \subfloat[$\eR=1.2$, $\ethe=20^\circ$, $N=500$.]%
    {\includegraphics[width=0.5\linewidth, trim = 0pt 0pt 0pt 100pt, clip]{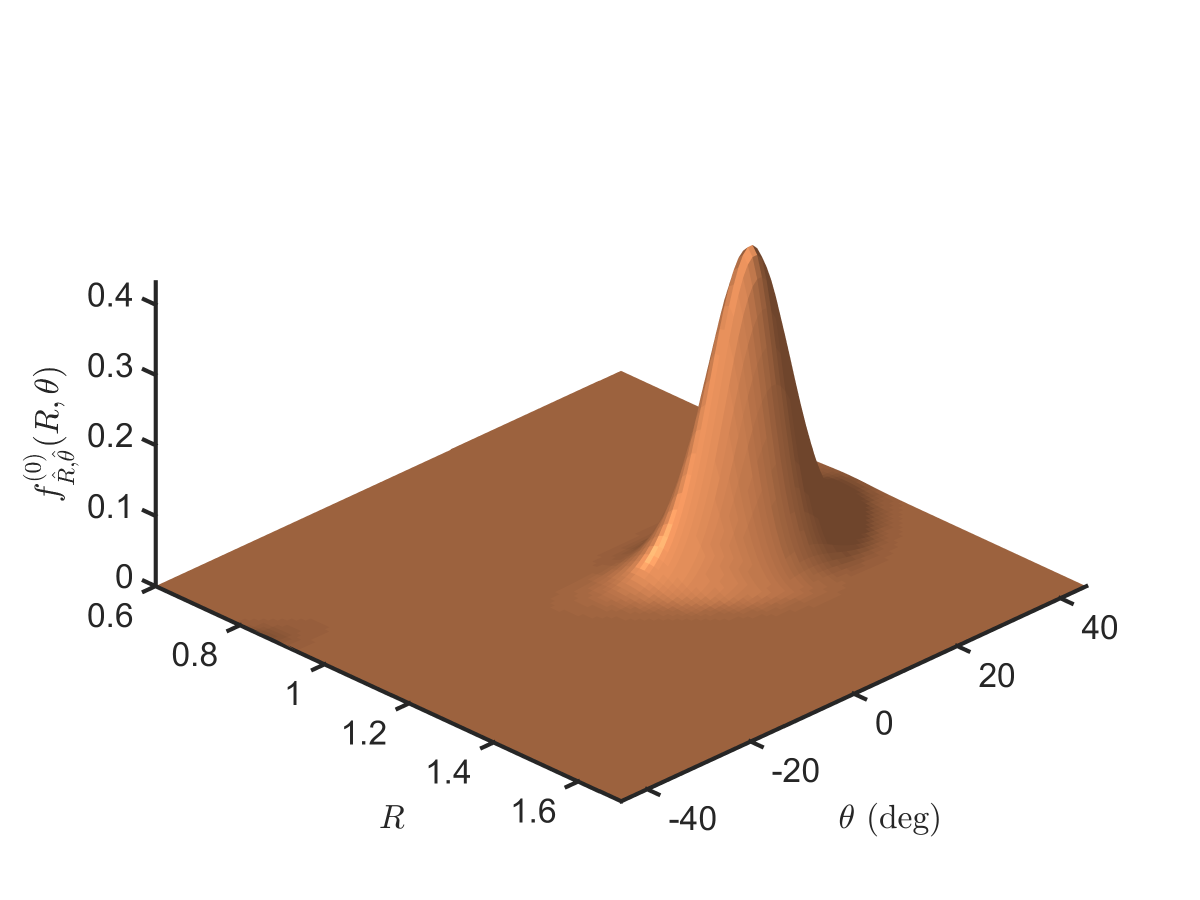}}\\
    \subfloat[$\eR=3$, $\ethe=10^\circ$, $N=100$.]%
    {\includegraphics[width=0.5\linewidth, trim = 0pt 0pt 0pt 100pt, clip]{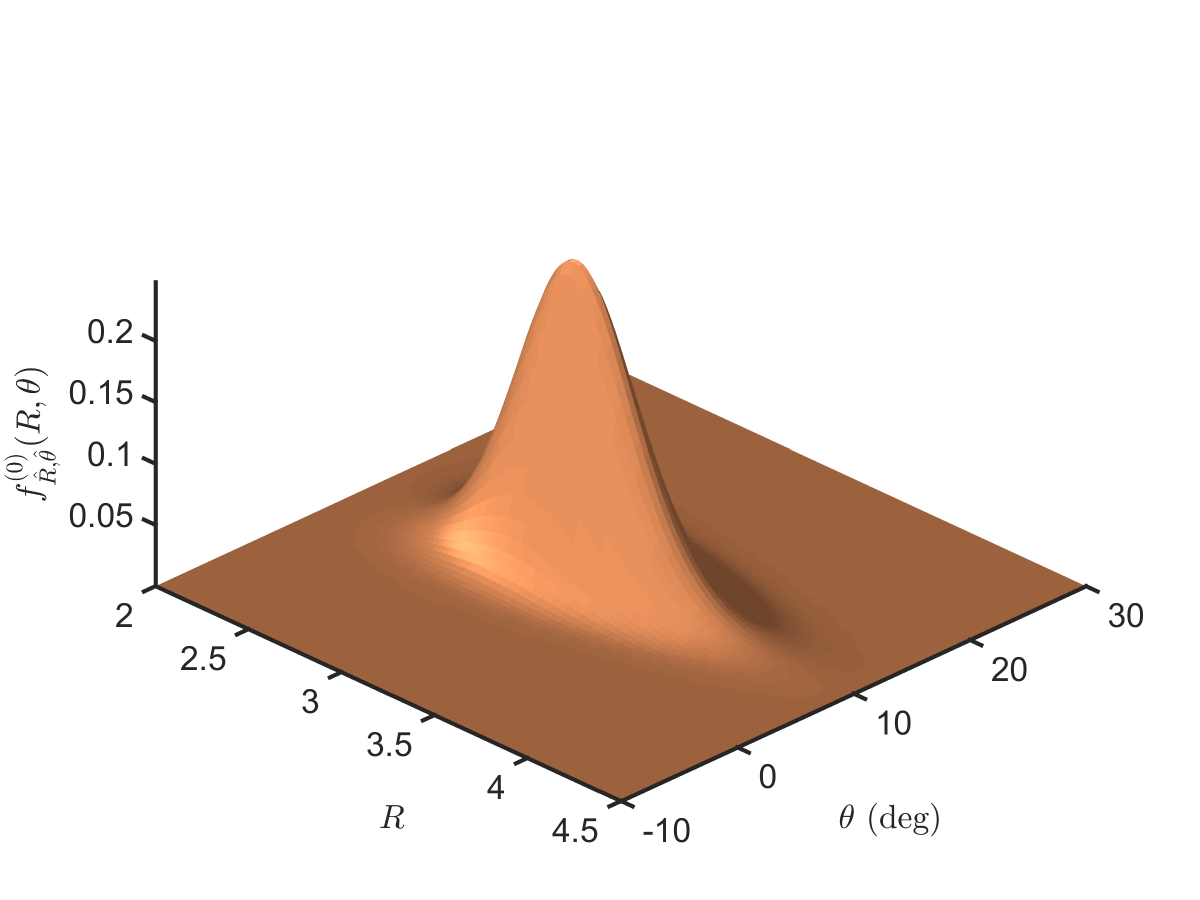}}
    \subfloat[$\eR=3$, $\ethe=10^\circ$, $N=500$.]%
    {\includegraphics[width=0.5\linewidth, trim = 0pt 0pt 0pt 100pt, clip]{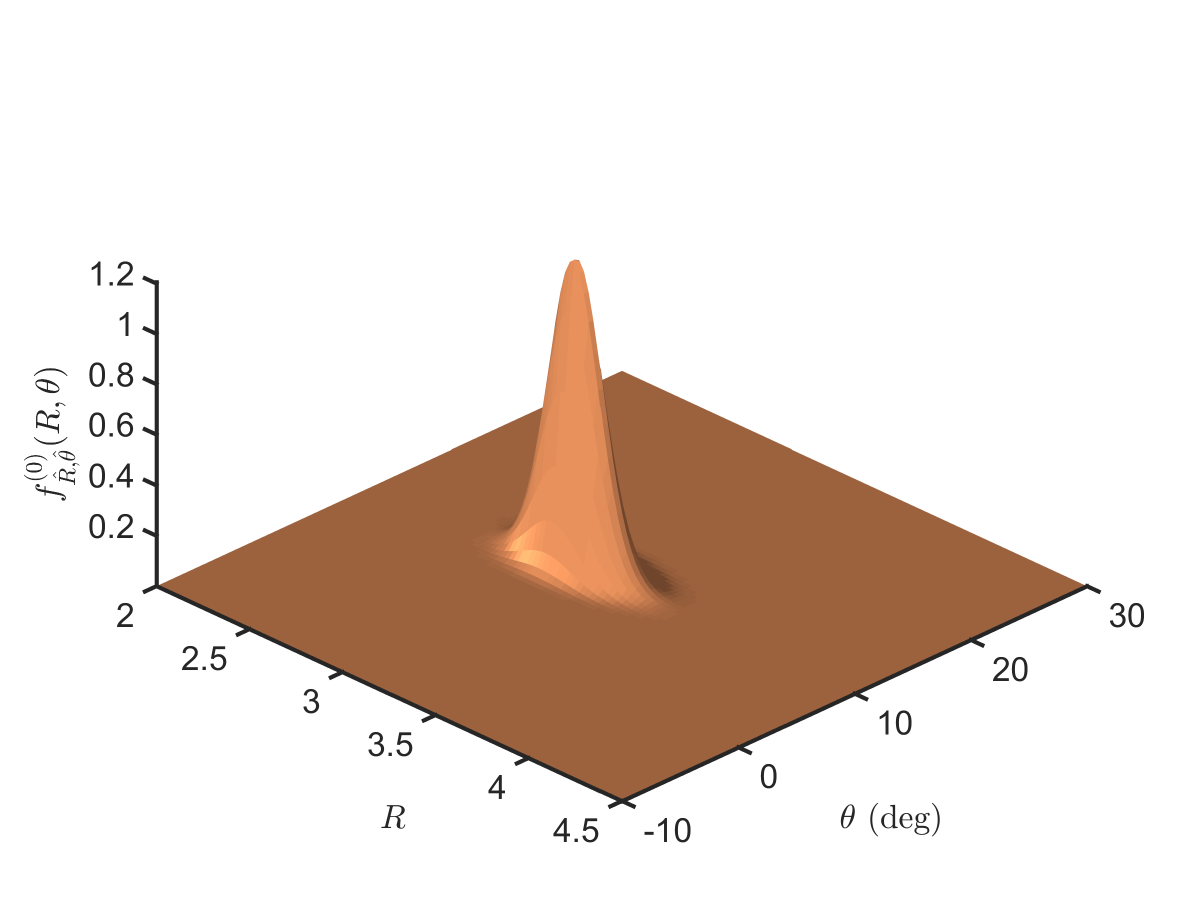}}
    \caption{Non-parametric JPDF $\NPJPDF$ for various anisotropy parameters $\eR,\ethe$ and sample size $N$.}
    \label{fig:JPDFs}
\end{figure}

Figure~\ref{fig:JPDFs} demonstrates representative plots of the non-parametric JPDF based on~\eqref{eq:f-anis_np_both}. Note the bimodal structure of the JPDF for $N=100$ in Fig.~\ref{subfig:JPDF_bimodal}, with one mode at $R=1.2$ and the other (smaller) at $R \approx 0.8$. This is due to the considerable spread of $\htheta$, which results from the relatively small number of sampling points and the degeneracy of the anisotropy vector, i.e., the fact that the combination $(R,\theta)$ is equivalent to $(1/R,\theta -\pi/2)$; the degenerate peak at $(0.83, -70^{\circ})$ is folded into the primary domain. On the other hand, the smaller dispersion of $\htheta$ for $R=3$ leads to a single mode even for $N=100$.

Knowledge of the anisotropy JPDF allows the construction of probability regions for the anisotropy statistics and confidence regions for the anisotropy parameters. For an SRF $X(\bs{s})$ that satisfies the conditions of Lemma~\ref{theor:fQ}, the probability region corresponding to level $p$ of the anisotropy statistics (based on the simplifications of the non-parametric approximation), is given by the following slight modification of~\eqref{eq:confi-reg}
\begin{equation}
    \label{eq:cr-R-theta}
    y_{0}^2(R, \theta; \eR, \ethe) - \frac{1}{2}  = \frac{\ln(1-p)}{N},
\end{equation}
where $y_{0}$ is a function of the values $R, \theta$, and the parameters $\eR, \ethe$ as defined in~\eqref{eq:y0}.


\section{Statistical Test of Isotropy}
\label{sec:isot-test}

\begin{theorem}[Isotropic ratio]
	\label{theor:isot-test-1}
	Let $X(\bfs)$ be a statistically isotropic GSRF $(\eR=1)$ which is sampled at $N$ points. Assume that the covariance $\cxx({\bf r})$ is short-ranged and its  spectral density satisfies $\tilde{C}({\bfk}) \sim \mathcal{O}(\kk^{-3-\epsilon})$ for $\epsilon>0$ as $\kk \ra \infty$ as defined in Lemma~\ref{theor:fQ}.
	In addition, assume that the asymptotic regime conditions hold. The probability interval of the anisotropic ratio at probability level $p$ (for $N > 2\ell_{p}$) is given by
	\begin{equation}
	    \label{eq:isot-R-int-np}
	     \left( \frac{N - 2\sqrt{\ell_p (N - \ell_p)}}{N - 2\ell_p},
	    \frac{N + 2\sqrt{\ell_p (N - \ell_p)}}{N - 2\ell_p} \right),
	\end{equation}
	
	\noindent where $\ell_{p}=F^{-1}(\chi^{2}=p,\nu=2)=-2\ln(1-p)$ is the inverse of the chi square cumulative distribution function with two degrees of freedom.
\end{theorem}

\begin{proof}
	For $\eR=1$  the angle dependent terms in the equations~\eqref{eqs:coeff} vanish, showing explicitly that the probability region is independent of $\theta$. Plugging~\eqref{eq:y0} in~\eqref{eq:cr-R-theta}  the following quadratic in $R^2$ equation is obtained
	\begin{equation*}
		 N\,(R^2-1)^2 - 2\ell_p \,(R^4+1)=0.
	\end{equation*}
	
	In fact, the probability region is reduced to a one-dimensional probability interval whose endpoints coincide with the roots of the above equation. The constraint $N>2\ell_{p}$ is in practice satisfied for $N \to \infty$	and ensures that the roots of the above equation  are positive real numbers. Based on the definition of $\ell_{p}$ the constraint is equivalent to $ N > -4\ln(1-p)$. For example, $p=0.95$  implies $\ell_p \approx 6$ and $N>12$.
\end{proof}

\begin{figure}
    \centering
    \includegraphics[width=0.66 \linewidth]{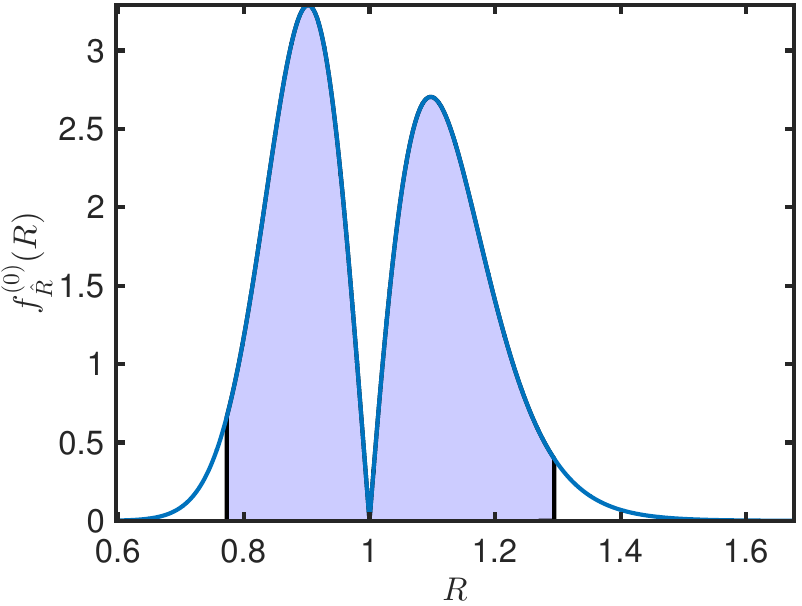}
    \caption{Non-parametric sampling PDF of the anisotropic ratio, $f_{\hR}^{(0)}(R)$, for an isotropic random field sampled at $N=100$  points. Shaded area represents the corresponding $95\%$ probability interval $(R_{-},R_{+}) = (0.77,1.29)$.}
    \label{fig:ci}
\end{figure}

Equation~\eqref{eq:isot-R-int-np} is independent of $\cxx(\bfr)$ and thus provides a non-parametric approximation of the probability interval for $\eR$. The JPDF~\eqref{eq:f-anis_np_both} is independent of $\theta$ and $\ethe$ for  $\eR=1$. The PDF, $f_{\hR}^{(0)}(R)$, of $\hR$ for $\eR=1$ and $N=100$ is shown in Fig.~\ref{fig:ci}, including  the $95\%$ probability interval predicted by~\eqref{eq:isot-R-int-np}.
Note that the PDF has a node instead of a peak at $R=1$. This is not an artifact of the non-parametric approximation, since the complete JPDF~\eqref{eq:f-anispar} also vanishes at $R=1$. The node is due to the root of the Jacobian~\eqref{eq:J} at $R=1$, which reflects that the isotropic point $(1,0)$ in $(\qd, \qo)$-space is mapped onto the straight line  $R=1$ in the $(R,\theta)$-space. The node is also evident in numerical simulations that do not use the Jacobian (see Figure~\ref{fig:scattered_iso} below).


\section{Application to Simulated and Real Data}
\label{sec:simulations}
To apply the formalism developed above to data sets that comprise discrete sets of values, we replace the partial derivatives by respective discrete operators $\check{\partial}{X}_{i}(\bs{s}_k)$,  $i=1,2$. The respective estimates of $\bs{\eQ}$ are denoted by $\bs{\cQ}$. The discretization introduces a bias that increases with the sparsity of the sampling pattern. A ``good'' sampling pattern is characterized by a typical distance $\hat{a}$ between nearest neighbors which is approximately uniform  (ideally, a regular lattice pattern is best) and $\hat{a} \ll \min(\xi_1,\xi_2)$, where $\xi_1,\xi_2$ are the principal correlation lengths. Different approaches for estimating $\check{\partial}{X}_{i}(\bs{s}_k)$ are investigated in~\citep{dth08}. Herein, the centered differences scheme is used for gridded data.

We denote average values of a statistic over different samples (repetitions) by a bar over the respective symbol, i.e., $\overline{\cQ}_{ij}$. For simulated data, the ensemble properties $\bs{\eQ}$ and $\CQQ$ which are unknown \emph{a priori}, are replaced by the respective averages $\bs{\eQ} \approx (\overline{\cQ}_{11},\overline{\cQ}_{22},\overline{\cQ}_{12})^{t}$ and ${\CQQ} \approx \bs{C}_{\bs{\overline{\cQ}}}$. In the non-parametric approximation, $\CQQ^{(0)}$ is obtained from~\eqref{eq:Cvv-nonpar} by replacing $\bs{\eQ}$ with $\overline{\bs{\cQ}}$.

\subsection{Simulated Scattered Data}
\label{ssec:simulations-B}

\begin{figure*}
	\captionsetup[subfigure]{margin = 5pt}    
	\centering
	\subfloat[Gaussian SRF, $\eR=1.5$, $\ethe=-30^\circ$]{\label{subfig:GL} \includegraphics[width=0.33\linewidth, trim = 35pt 0pt 75pt 20pt, clip]{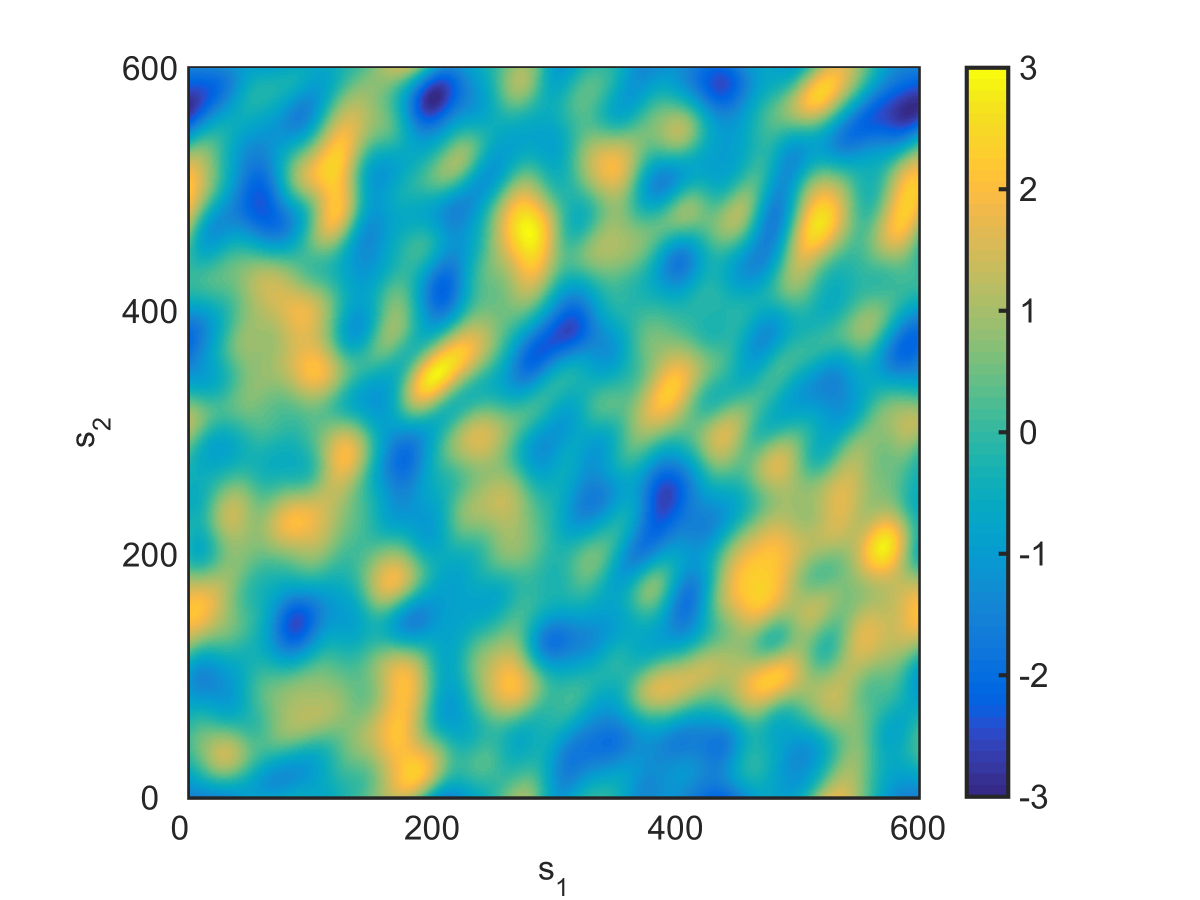}} %
	\subfloat[1296 random nodes and interpolated field.]{\label{subfig:GS} \includegraphics[width=0.33\linewidth, trim = 35pt 0pt 75pt 20pt, clip]{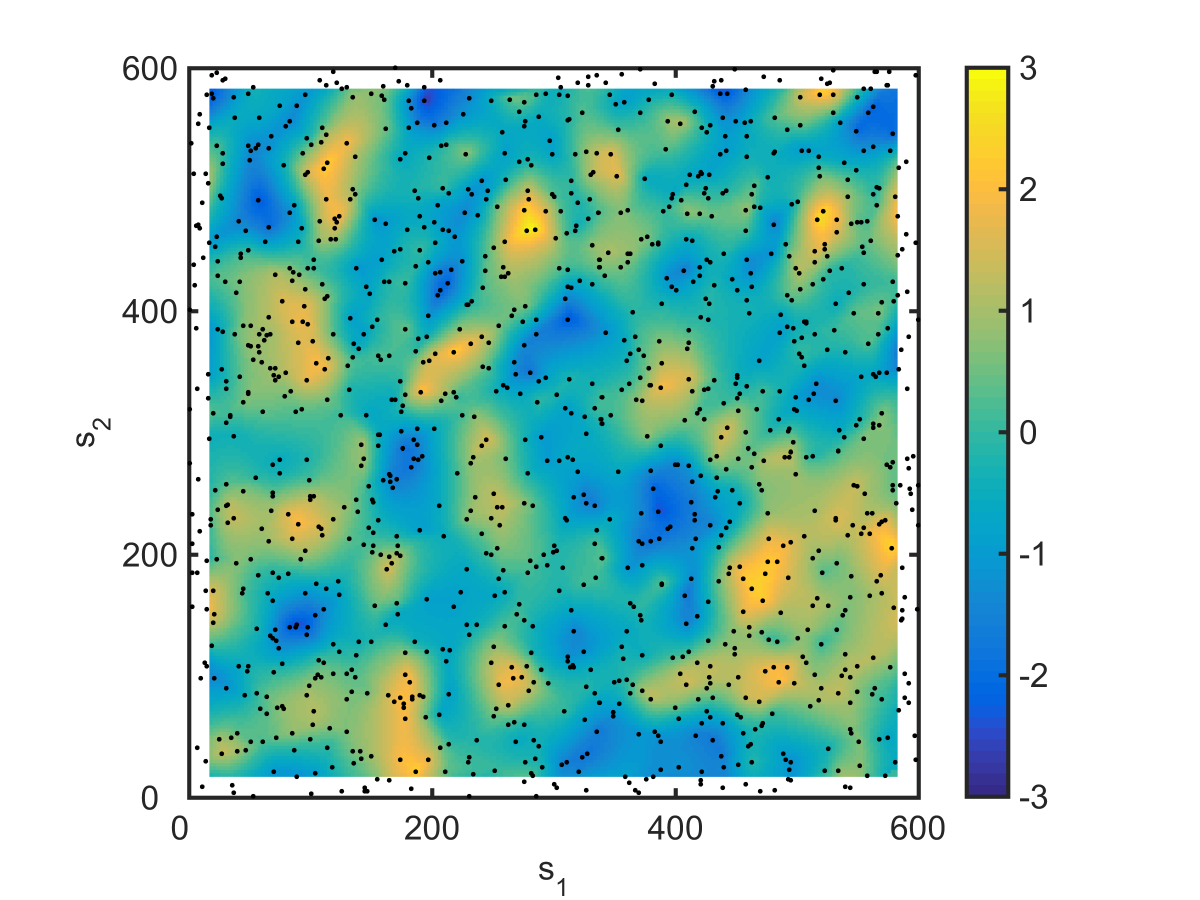}} %
	\subfloat[Anisotropy estimates and probability region]{\label{subfig:GSC}	\includegraphics[width=0.33\linewidth, trim = 35pt 0pt 75pt 0pt, clip]{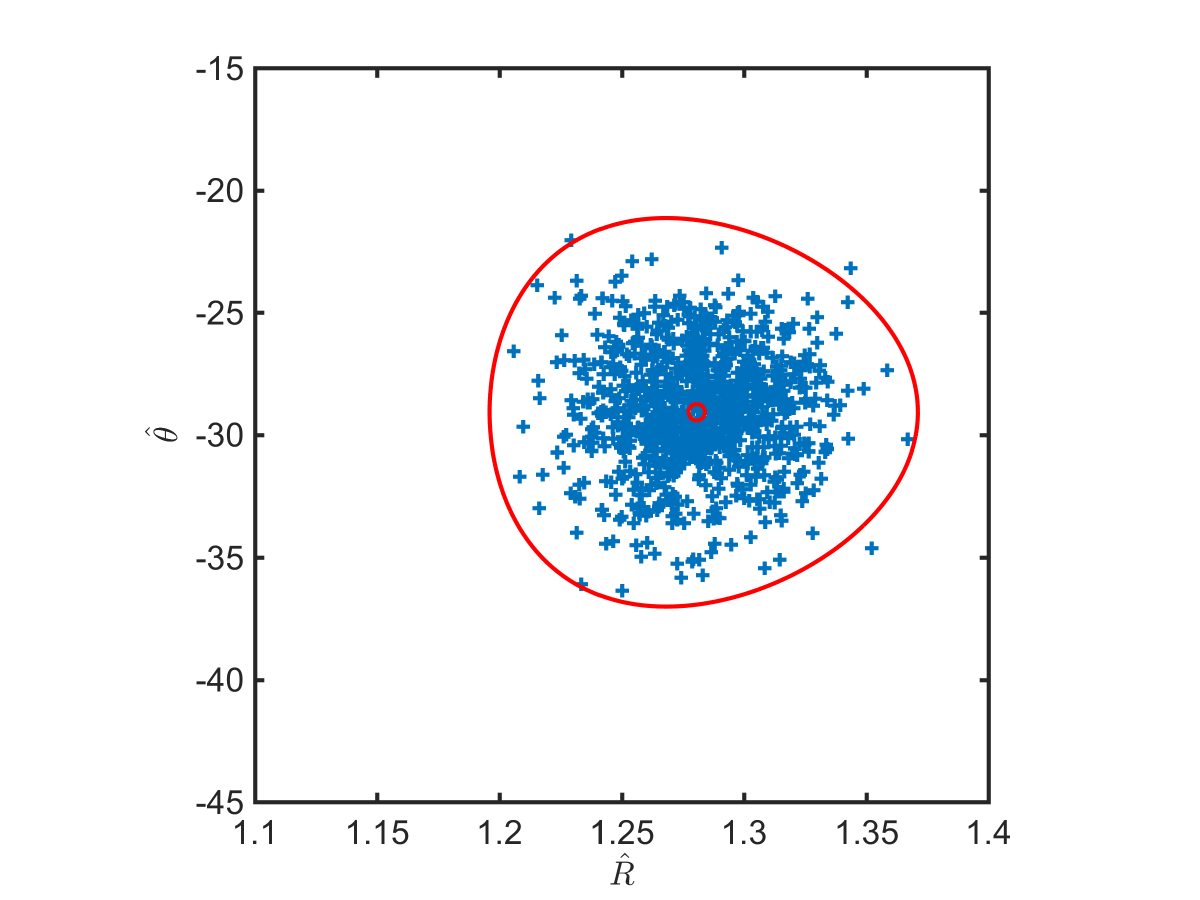}} %
	\\
	\subfloat[Mat\'ern SRF, $\eR=1.5$, $\ethe=-30^\circ$]{\label{subfig:ML} \includegraphics[width=0.33\linewidth, trim = 35pt 0pt 75pt 20pt, clip]{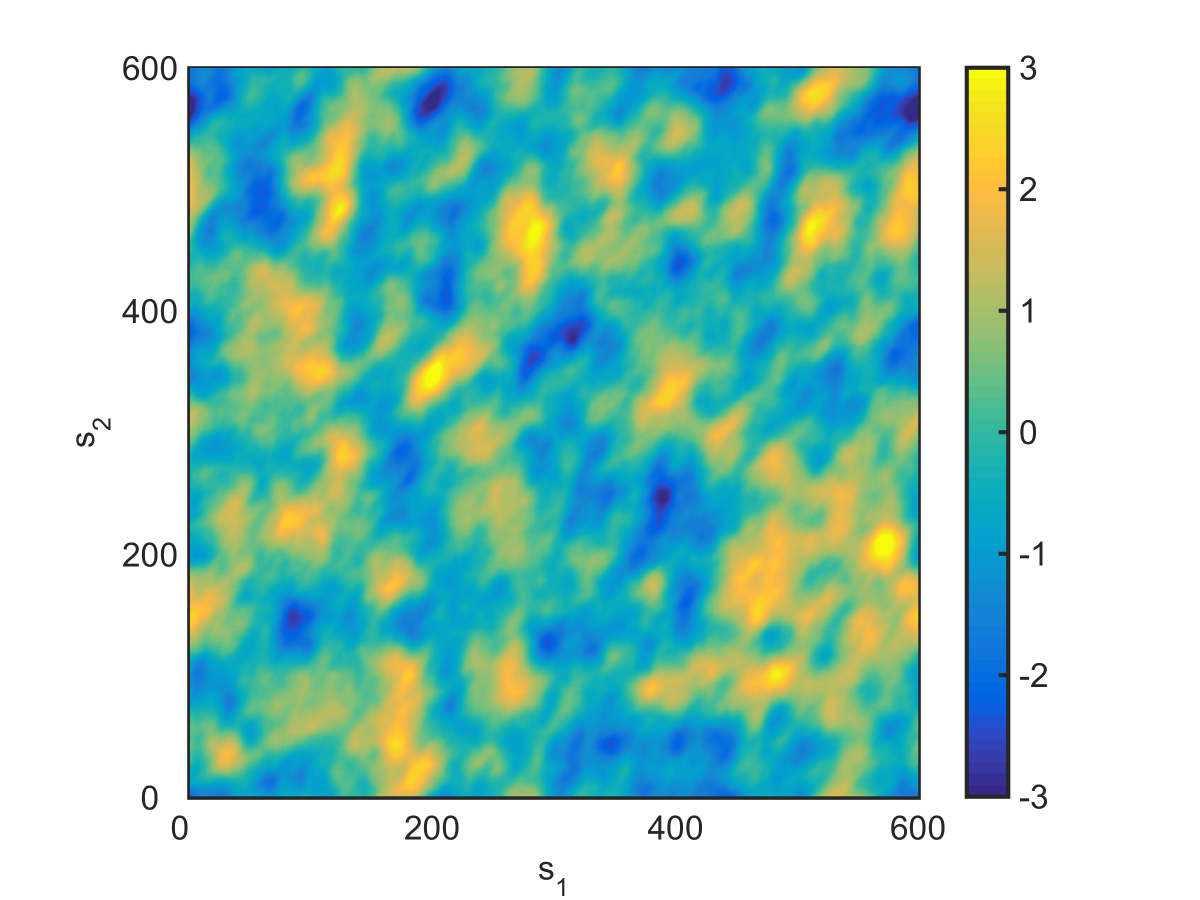}} %
	\subfloat[1296 random nodes and interpolated field.]{\label{subfig:MS} \includegraphics[width=0.33\linewidth, trim = 35pt 0pt 75pt 20pt, clip]{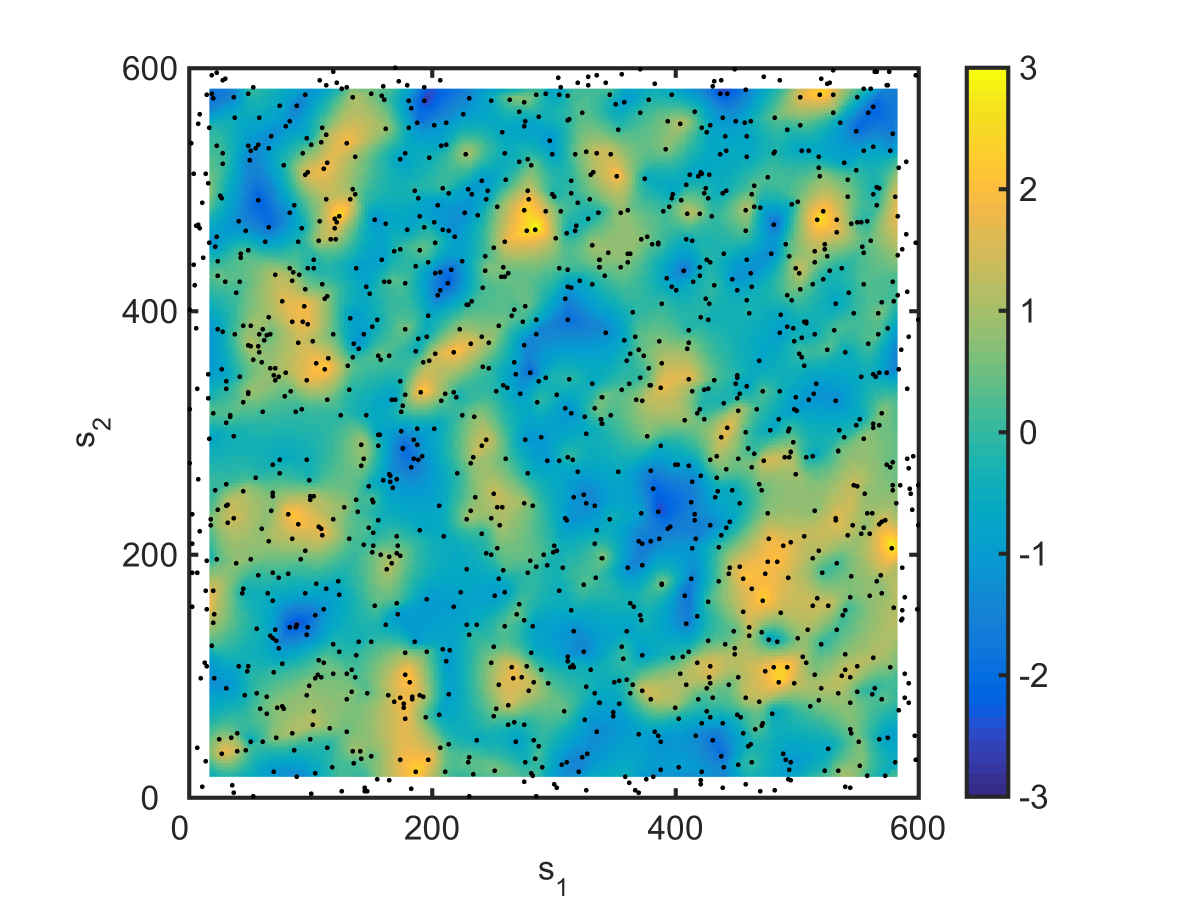}} %
	\subfloat[Anisotropy estimates and probability region]{\label{subfig:MSC}	\includegraphics[width=0.33\linewidth, trim = 35pt 0pt 75pt 0pt, clip]{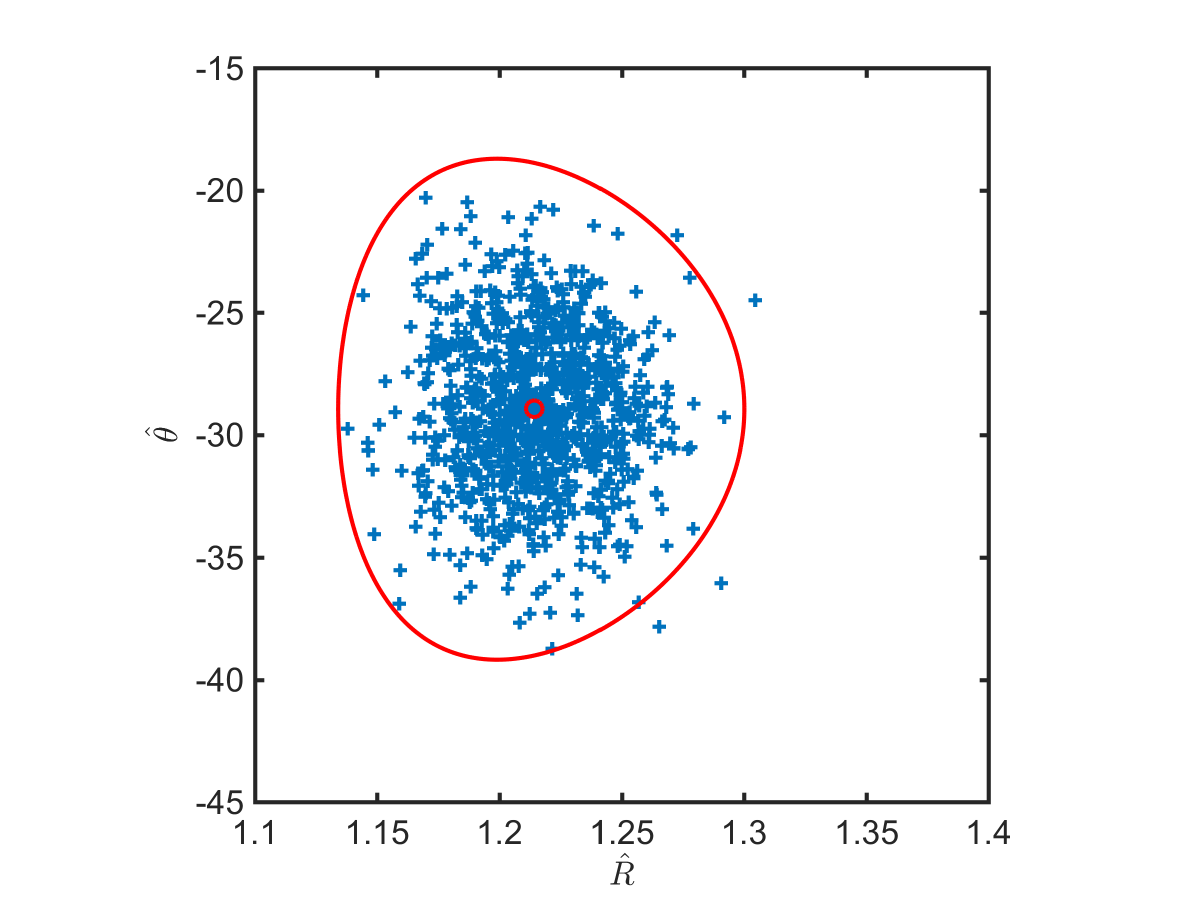}} %
	\caption{Non-parametric probability region estimation for scattered data.
	(a) and (d):~Realization of zero-mean, unit-variance anisotropic Gaussian SRF with $\xi=28.3$ and M\'atern SRF with $\nu=2$ and $\xi=10$ on a $600 \times 600$ square grid. (b) and (e):~Random sample of $N=1296$ points and interpolated field on a $200 \times 200$ grid using natural neighbors. Boundary strips of thickness $L/\sqrt{N}$ were discarded from the interpolated field  to avoid interpolation artifacts at the domain boundary. (c) and (f):~Anisotropy estimates (crosses) are generated from $1000$ random samples of $N=1296$ points; the continuous curve corresponds to  $95\%$ non-parametric probability region calculated with the ensemble-based anisotropy estimates $(\overline{\cR},\ \overline{\ctheta})$, which are denoted by a small circle inside the cloud.}
	\label{fig:scattered_aniso}
\end{figure*}

We generate SRF realizations with specified $(\eR,\ethe)$ to validate the probability region of the anisotropy parameters~\eqref{eq:cr-R-theta}. Figure~\ref{fig:scattered_aniso} and Table~\ref{t:validation_ensemble} investigate the anisotropic case $\eR=1.5$, $\ethe=-30^\circ$, whereas the isotropic case is considered in Table~\ref{t:validation} and Figure~\ref{fig:scattered_iso}. A desktop computer with an Intel\textsuperscript{\textregistered} Core\textsuperscript{\texttrademark} i5-2500 (4 cores, 3.30 GHz) CPU running \matlab\ R2015b under 64-bit Windows\textsuperscript{\textregistered} 7, was used for all the simulations.

We simulate scattered data using the following method: First, a realization of an GSRF is generated on a regular grid. The Fourier Filtering Method~\citep{Padro93, Lantu02, dth05} is used on $L \times L$ square grids with lattice constant $a=1$. We use Gaussian, $\cxx(\bs{r})=\stdx^{2} \, \exp(-\rr^2/\xi^2)$, and Mat\'ern, $\cxx(\bs{r}) = \stdx^{2} \, 2^{1-\nu} \, \Gamma(\nu)^{-1} \, \xi^{-\nu} \times \allowbreak \rr^{\nu} K_{\nu}(\rr/\xi)$, covariance functions (expressions correspond to the isotropic case), where $\Gamma(\cdot)$ is the Gamma and $K_{\nu}(\cdot)$ the modified Bessel function of order $\nu$. In the Gaussian case, the correlation range is controlled by $\xi$ whereas in the Mat\'ern case by both $\xi$ and $\nu$. The smoothness parameter $\nu$ adjusts the differentiability of the SRF: $\nu=1/2$  corresponds to the non-differentiable exponential function and $\nu \to \infty$ to the infinitely differentiable Gaussian. For given $\xi$, the field is smoother for higher $\nu$. To compensate for this effect and to compare SRFs of similar spatial variability, we use rescaled correlation lengths $\tilde{\xi} = A_d \xi$, where $A_d$ is the \emph{integral scale factor}~\citep{dth11:SERRA}: In $d=2$, $A_d = 2\sqrt{\pi \nu}$ for Mat\'ern correlations whereas for Gaussian correlations $A_d = \sqrt{\pi}$. For equal rescaled correlation lengths, $\tilde{\xi}_\text{Gauss} = \tilde{\xi}_\text{Mat\'ern}$, with $d=2$, $\nu=2$, it follows that $\xi_\text{Gauss} = 2\sqrt{2}\,\xi_\text{Mat\'ern}$.

We randomly choose a fraction of the grid points to mimic scattered data. For a square lattice of side $L$ a sample of $N = (\rho L)^2$ points are randomly chosen from Gaussian and Mat\'ern lattice SRFs. An estimate of the mean distance between $N$ uniformly distributed points is $L/\sqrt{N} = 1/\rho$, thus $\rho$ is the mean sampling frequency. The samples respect the condition that the correlation lengths exceed the mean distance between the points, as specified in the first paragraph of this Section.

We employ the natural neighbor interpolation method~\citep{Ledoux05} in \matlab\ on an $M \times M$ square grid with $M=200$. Natural neighbor interpolation provides smooth surfaces and does not assume isotropy of the data; however, it is defined only inside the convex hull of the data sites. Due to the occasionally poor sampling near the domain boundaries, interpolation artifacts appear~\citep{Bobach2009} as elongated islands, oriented vertically along the left and horizontally along the bottom sides of the domain. Hence, they tend to bias the anisotropy estimates towards higher or lower anisotropy ratios and angles near zero. Thus, boundary strips of thickness $L/\sqrt{N}$ are discarded from the interpolation surface to minimize bias. The partial derivatives are estimated via centered differences on the interpolated surface. 
Finally, we perform anisotropy estimation for each sample and compute the non-parametric probability region at $p=0.95$ using ensemble averages. Also, we compute confidence regions for each anisotropy estimate at several confidence levels.

\subsubsection{Anisotropic Scattered Data}
Figure~\ref{subfig:GL} demonstrates a realization of a zero-mean, unit variance anisotropic GSRF with Gaussian covariance with $\eR=1.5, \ethe=-30^\circ$, $\xi=28.3$ on a $600 \times 600$ grid. A randomly extracted set of $N=1296$ points ($\rho = 0.06$) is shown in Fig.~\ref{subfig:GS}. The depicted smooth field is generated from the 1296 points by interpolation and is used to estimate $(\cR, \ctheta)$.

In Figure~\ref{subfig:GSC} the non-parametric probability region (red contour) at $p=0.95$, defined by~\eqref{eq:cr-R-theta} is compared with CHI anisotropy estimates (blue crosses) from $1000$ SRF samples. For each sample, we estimate $\bs{\eQ}$ by means of the spatial average $\bs{\cQ}$ and then calculate $(\cR, \ctheta)$ by applying~\eqref{eq:R-theta}. We estimate $(\eR,\ethe)$ based on the $\overline{\cR},\ \overline{\ctheta}$, which are obtained from the ensemble average $\overline{\bs{\cQ}}$ by means of Theorem~\ref{theor:aniso}. The ensemble-based anisotropy estimate $(\overline{\cR},\ \overline{\ctheta})$, is denoted by a small circle inside the cloud of the $(\cR, \ctheta)$ points. Figures~\ref{subfig:ML}--\ref{subfig:MSC} demonstrate the simulated scattered data probability region estimation for a zero-mean, unit-variance Mat\'ern covariance with $\eR=1.5, \ethe=-30^\circ$, $\xi = 10$, $\nu=2$. The normality of $\Qind$, supported by CLT considerations as shown in Lemma~\ref{theor:fQ}, was confirmed by normal probability plots (not shown here).

The non-parametric probability region (Theorem~\ref{theor:JPDF-NP}) extends beyond the region obtained from the true JPDF (this is supported by Figures \ref{subfig:GSC} and \ref{subfig:MSC} as explained below). We conducted numerical experiments (not shown here) for several values of $\tilde{\xi}/a$ and $N$ to confirm that non-parametric probability regions based on~\eqref{eq:f-anis-np-jac} are more extended in parameter space than the regions based on the true JPDF~\eqref{eq:triv-pdf-scaled}. If $\tilde{\xi}/a \to 0$, i.e., as the spatial extent of the correlations is reduced, the scatter cloud of $(\cR, \ctheta)$ expands and tends to fill the non-parametric probability region. On the other hand, as $\tilde{\xi}/a $ increases, i.e., for dense sampling of the SRF, the scatter cloud tends to be confined inside the smaller parametric  region. These observations agree with our earlier statement that the non-parametric approximation contains the true probability region.

In Table~\ref{t:validation_ensemble} we validate the non-parametric anisotropy \emph{confidence region} for simulated scattered Gaussian ($\xi=28.3$) and Mat\'ern ($\xi=10$, $\nu=2$) covariance functions with $\eR=1.5$ and $\ethe=-30^\circ$. We generate 1000 realizations for different domain sizes ($L=600$, $800$, $1000$, $1200$) and mean sampling frequencies ($\rho=0.04, 0.06$) and we enumerate the number of simulations for which the ensemble means $\overline{\check{R}},\overline{\check{\theta}}$ (as estimates of the population means) are outside the non-parametric confidence region. The latter is computed for each anisotropy estimate at different confidence levels ($p=0.95$, $0.75$, $0.68$, $0.5$, $0.25$) using~\eqref{eq:cr-R-theta}. If the true JPDF and the confidence regions of the anisotropy statistics are known at the $p$ levels above, the average number of simulations for which the true confidence region does not contain the ensemble means $\overline{\check{R}},\overline{\check{\theta}}$ is $50$, $250$, $320$, $500$, and $750$ respectively. However, the number of simulations for which $\overline{\check{R}},\overline{\check{\theta}}$ lie outside the non-parametric region~\eqref{eq:cr-R-theta} is always less than expected for the true confidence regions. This observation agrees with the proposition that the non-parametric confidence region~\eqref{eq:cr-R-theta} contains the true confidence region.

\begin{sidewaystable}
        \renewcommand{\arraystretch}{1.3}
        \setlength{\tabcolsep}{3pt}
        \small
        \centering
        \caption{Numerical validation of the non-parametric anisotropy confidence region for simulated scattered data. Anisotropy estimates were computed from 1000 random samples of $N$ points from (a) Gaussian and (b) Mat\'ern lattice SRFs with $R=1.5$, $\theta=-30^\circ$. For each lattice size $L$ and sampling frequency $\rho$, $N_{\text{out},p}$ is the number of anisotropy estimates for which the non-parametric confidence region computed at $p=0.95$, $0.75$, $0.68$, $0.5$, and $0.25$ using Theorem~\ref{theor:JPDF-NP}, does not contain the ensemble mean $(\overline{\check{R}},\overline{\check{\theta}})$.}
        \label{t:validation_ensemble}
        \subfloat[Gaussian,  $R=1.5$, $\theta=-30^\circ$, $\xi = 28.3$, $\rho=0.04,0. 06$]{%
                \label{t:gaussian_ensemble} %
                \begin{tabular}{lcccccccc}
                	\hline\hline
                	$L$                                              &        \multicolumn{2}{c}{$600$}         &        \multicolumn{2}{c}{$800$}         &         \multicolumn{2}{c}{$1000$}         &        \multicolumn{2}{c}{$1200$}        \\
                	\hline
                	$N$                                              &       $576$        &       $1296$        &       $1024$        &       $2304$       &        $1600$        &       $3600$        &       $2304$       &       $5184$        \\
                	$N_{\text{out},p=0.95}$                          &        $12$        &         $3$         &        $10$         &        $3$         &         $1$          &         $0$         &        $6$         &         $1$         \\
                	$N_{\text{out}, p=0.75}$                         &       $132$        &        $68$         &        $106$        &        $41$        &         $83$         &        $26$         &        $85$        &        $20$         \\
                	$N_{\text{out}, p=0.68}$                         &       $194$        &        $102$        &        $159$        &        $68$        &        $141$         &        $48$         &       $128$        &        $47$         \\
                	$N_{\text{out}, p=0.50}$                         &       $380$        &        $280$        &        $324$        &       $191$        &        $288$         &        $155$        &       $274$        &        $166$        \\
                	$N_{\text{out}, p=0.25}$                         &       $658$        &        $560$        &        $601$        &       $534$        &        $592$         &        $465$        &       $574$        &        $450$        \\
                	$\overline{\check{R}},\overline{\check{\theta}}$ & $1.19,-29.0^\circ$ & $1.28, -29.8^\circ$ & $1.21, -29.8^\circ$ & $1.30,-29.7^\circ$ & $1.22,  -30.1^\circ$ & $1.31, -30.1^\circ$ & $1.22,-30.4^\circ$ & $1.32, -30.9^\circ$ \\
                	\hline\hline
                \end{tabular}%
        }\\
        \setlength{\tabcolsep}{3pt}
        \subfloat[Mat\'ern, $R=1.5$, $\theta=-30^\circ$, $\nu=2$, $\xi=10$,$\rho=0.04, 0.06$]{%
                \label{t:matern_ensemble} %
                \begin{tabular}{lcccccccc}
                	\hline\hline
                	$L$                                              &        \multicolumn{2}{c}{$600$}         &        \multicolumn{2}{c}{$800$}         &         \multicolumn{2}{c}{$1000$}         &        \multicolumn{2}{c}{$1200$}        \\
                	\hline
                	$N$                                              &       $576$        &       $1296$        &       $1024$        &       $2304$       &        $1600$        &       $3600$        &       $2304$       &       $5184$        \\
                	$N_{\text{out},p=0.95}$                          &        $8$         &         $4$         &        $10$         &        $3$         &         $5$          &         $2$         &        $11$        &         $4$         \\
                	$N_{\text{out}, p=0.75}$                         &       $150$        &        $78$         &        $130$        &        $69$        &         $95$         &        $42$         &       $114$        &        $51$         \\
                	$N_{\text{out}, p=0.68}$                         &       $209$        &        $121$        &        $172$        &       $120$        &        $138$         &        $78$         &       $170$        &        $86$         \\
                	$N_{\text{out}, p=0.50}$                         &       $391$        &        $289$        &        $346$        &       $257$        &        $302$         &        $206$        &       $320$        &        $227$        \\
                	$N_{\text{out}, p=0.25}$                         &       $680$        &        $608$        &        $667$        &       $572$        &        $632$         &        $515$        &       $630$        &        $543$        \\
                	$\overline{\check{R}},\overline{\check{\theta}}$ & $1.14,-28.1^\circ$ & $1.22, -28.9^\circ$ & $1.16, -29.5^\circ$ & $1.23,-29.8^\circ$ & $1.16,  -30.0^\circ$ & $1.24, -30.3^\circ$ & $1.17,-30.4^\circ$ & $1.25, -30.7^\circ$ \\
                	\hline\hline
                \end{tabular}%
        }
\end{sidewaystable}

\subsubsection{Isotropic Scattered Data}
We numerically validate the isotropy testing procedure by enumerating the number of anisotropy estimates that fall outside the probability region at $p=0.95$ for 1000 realizations of simulated scattered data in different domain sizes ($L=600,800,1000,1200$), mean sampling frequency ($\rho = 0.04, 0.06$), and isotropic covariances (Gaussian with  $\xi=28.3$, Mat\'ern with $\xi=10,\nu=2$). If the true JPDF and the probability regions of the anisotropy statistics are known at $p=0.95$ probability level, on average 50 out of the 1000 simulations should fall outside the true  region. In Table~\ref{t:validation}, $N_{\text{out,iso}}$ is the number of estimates that fall outside the $p=0.95$ isotropy probability interval $(R_{-},R_{+})$ using Eq.~\eqref{eq:isot-R-int-np} while $N_{\text{out}}$ is the number of  samples that fall outside the probability region calculated using the ensemble-based anisotropy estimate $(\overline{\cR},\ \overline{\ctheta})$ and Theorem~\ref{theor:JPDF-NP}. The mean time $\bar{t}$ for anisotropy estimation per processor core is also shown with an error estimate of one standard deviation.

\begin{sidewaystable}
	\renewcommand{\arraystretch}{1.3}
	\setlength{\tabcolsep}{3.55pt}
	\small
	\centering
	\caption{Numerical validation of the non-parametric isotropy test at $p=0.95$ for simulated scattered data. Anisotropy estimates generated from 1000 random samples of $N$ points from isotropic (a) Gaussian and (b) Mat\'ern lattice SRFs. For each lattice size $L$ and sampling frequency $\rho$, $N_{\text{out,iso}}$ is the number of simulations that fall outside the non-parametric probability region estimated from $\overline{\bs{\cQ}}$ using Theorem~\ref{theor:JPDF-NP}. $N_\text{out}$ is the number of simulations that fall outside the non-parametric isotropy  interval $(R_{-},R_{+})$ based on~\eqref{eq:isot-R-int-np}. $\bar{t}$: average time per anisotropy estimation (per processor core) with an error estimate of one standard deviation.}
	\label{t:validation}
	\subfloat[Gaussian,  $\xi = 28.3$, $\rho=0.04,0. 06$]{%
		\label{t:gaussian} %
		\begin{tabular}{lcccccccc}
			\hline\hline
			$L$                                              &    \multicolumn{2}{c}{$600$}    &    \multicolumn{2}{c}{$800$}    &    \multicolumn{2}{c}{$1000$}    &   \multicolumn{2}{c}{$1200$}    \\
			\hline
			$N$                                              &     $576$      &     $1296$     &     $1024$     &     $2304$     &     $1600$      &     $3600$     &     $2304$     &     $5184$     \\
			$N_{\text{out,iso}}$                             &      $14$      &      $4$       &      $12$      &      $2$       &       $6$       &      $0$       &      $11$      &      $1$       \\
			$N_{\text{out}}$                                 &      $12$      &      $3$       &      $13$      &      $2$       &       $5$       &      $0$       &      $9$       &      $1$       \\
			$\overline{\check{R}},\overline{\check{\theta}}$ & $0.990, -42.2$ & $0.986, -43.2$ & $0.995, -41.6$ & $0.992, -42.2$ & $0.997,  -43.0$ &  $1.00, 41.8$  & $0.998, -42.5$ &  $1.00, 44.9$  \\
			$R_{-}, R_{+}$                                   & $0.902, 1.11$  & $0.934, 1.07$  & $0.926,  1.08$ & $0.950, 1.05$  &  $0.940, 1.06$  & $0.960, 1.04$  & $0.950, 1.05$  & $0.967, 1.03$  \\
			$\bar{t}$ (msec)                                 &  $129 \pm 4$   &  $131 \pm 3$   &  $132 \pm 3$   &  $146 \pm 5$   &   $141 \pm 5$   &  $160 \pm 5$   &  $150 \pm 9$   &  $175 \pm 6$   \\
			\hline\hline
		\end{tabular}%
	}\\
	\setlength{\tabcolsep}{5pt} 
	\subfloat[Mat\'ern, $\nu=2$, $\xi=10$, $\rho=0.04, 0.06$]{%
		\label{t:matern} %
		\begin{tabular}{lcccccccc}
			\hline\hline
			$L$                                              &    \multicolumn{2}{c}{$600$}    &    \multicolumn{2}{c}{$800$}    &   \multicolumn{2}{c}{$1000$}    &   \multicolumn{2}{c}{$1200$}    \\
			\hline
			$N$                                              &     $576$      &     $1296$     &     $1024$     &     $2304$     &     $1600$     &     $3600$     &     $2304$     &     $5184$     \\
			$N_{\text{out,iso}}$                             &      $15$      &      $8$       &      $12$      &      $8$       &      $10$      &      $5$       &      $10$      &      $7$       \\
			$N_{\text{out}}$                                 &      $11$      &      $7$       &      $12$      &      $5$       &      $10$      &      $2$       &      $12$      &      $4$       \\
			$\overline{\check{R}},\overline{\check{\theta}}$ & $0.990, -44.4$ &  $1.01, 43.6$  &  $1.01,42.7$   & $0.993, -42.3$ &  $1.00, 43.6$  &  $1.01,44.2$   & $0.997, -37.6$ &  $1.00, 42.8$  \\
			$R_{-}, R_{+}$                                   & $0.902, 1.11$  & $0.934, 1.07$  & $0.926, 1.08$  & $0.950, 1.05$  & $0.940,  1.06$ & $0.960, 1.04$  & $0.950, 1.05$  & $0.967, 1.03$  \\
			$\bar{t}$ (msec)                                 &  $129 \pm 3$   &  $131 \pm 4$   &  $130 \pm 4$   &  $143 \pm 5$   &  $139 \pm 6$   &  $156 \pm 5$   &  $145 \pm 5$   &  $169 \pm 7$   \\
			\hline\hline
		\end{tabular}%
	}
\end{sidewaystable}

\begin{figure*}
	\captionsetup[subfigure]{margin = 15pt}    
	\centering
	\subfloat[Gaussian, $L=600$, $N=576$]{\label{subfig:G1} \includegraphics[width=0.33\linewidth, trim = 35pt 0pt 75pt 0pt, clip]{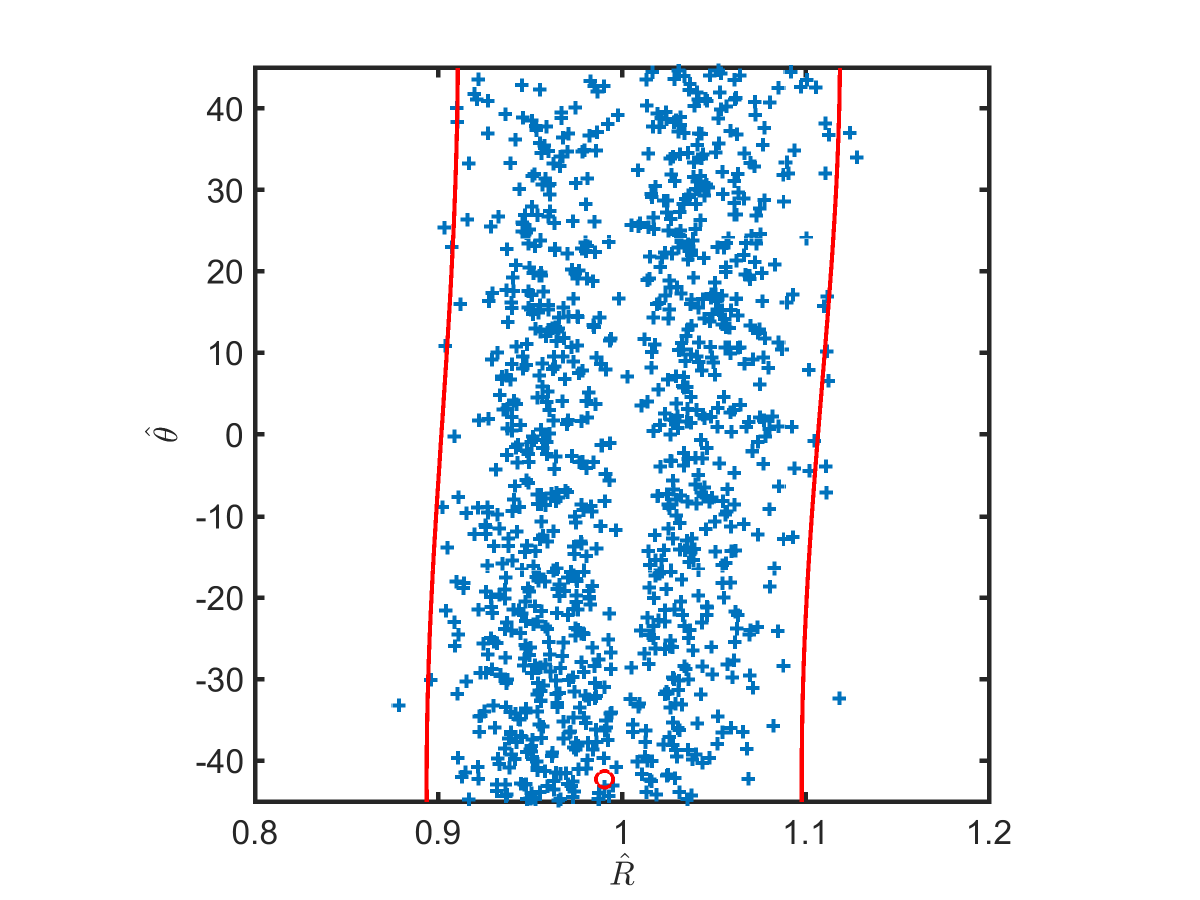}} %
	\subfloat[Gaussian, $L=800$, $N=1024$]{\label{subfig:G2} \includegraphics[width=0.33\linewidth, trim = 35pt 0pt 75pt 0pt, clip, clip]{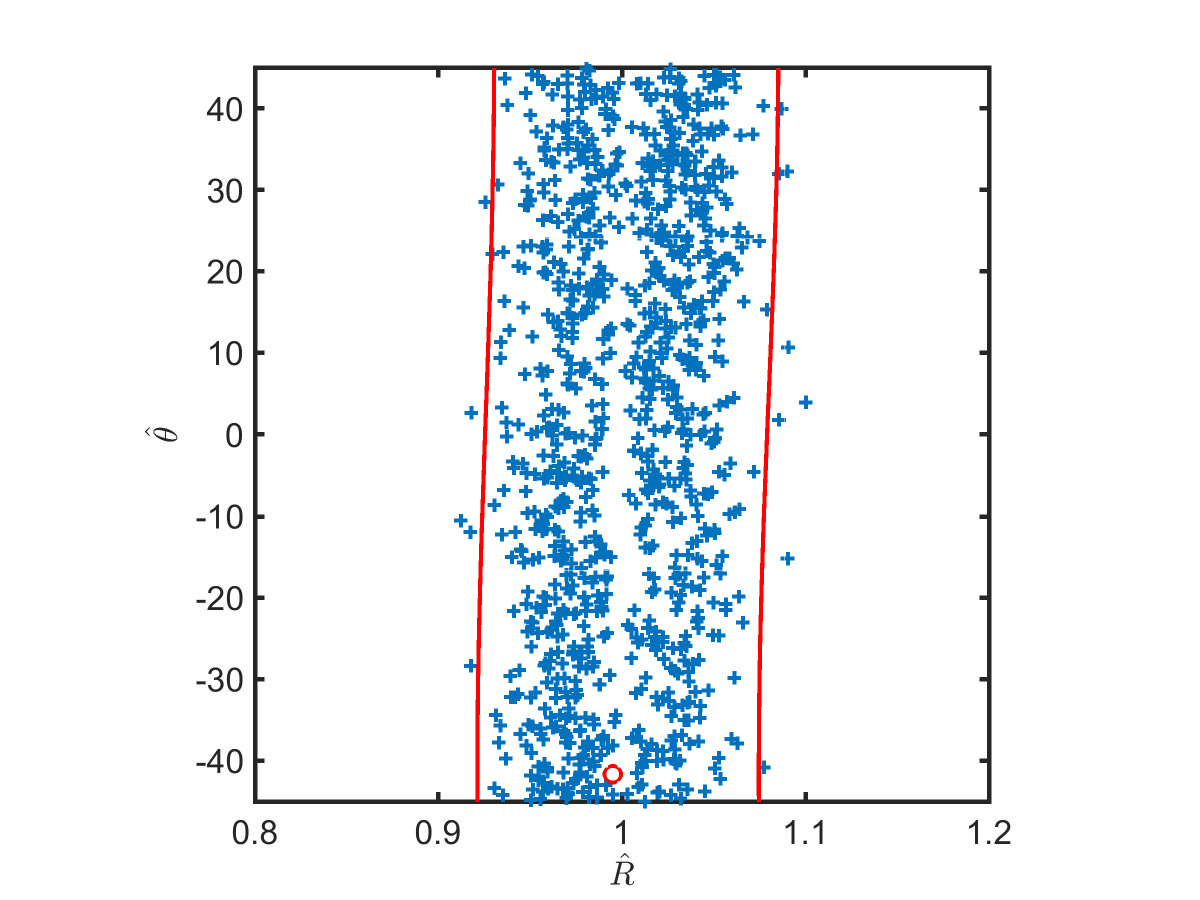}} %
	\subfloat[Gaussian, $L=1000$, $N=1600$]{\label{subfig:G3} \includegraphics[width=0.33\linewidth, trim = 35pt 0pt 75pt 0pt, clip, clip]{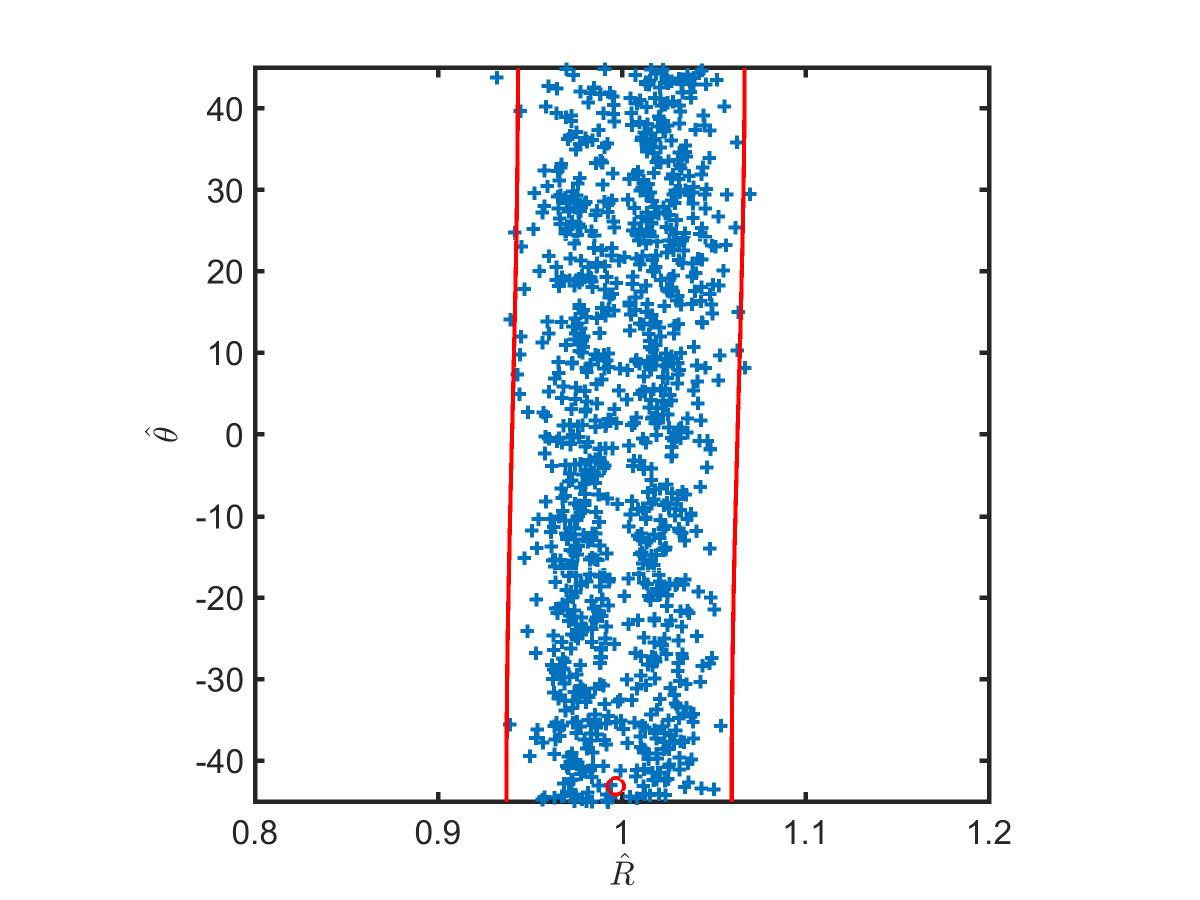}} %
	\\
	\subfloat[Mat\'ern, $L=600$, $N=576$]{\label{subfig:M1} \includegraphics[width=0.33\linewidth, trim = 35pt 0pt 75pt 0pt, clip]{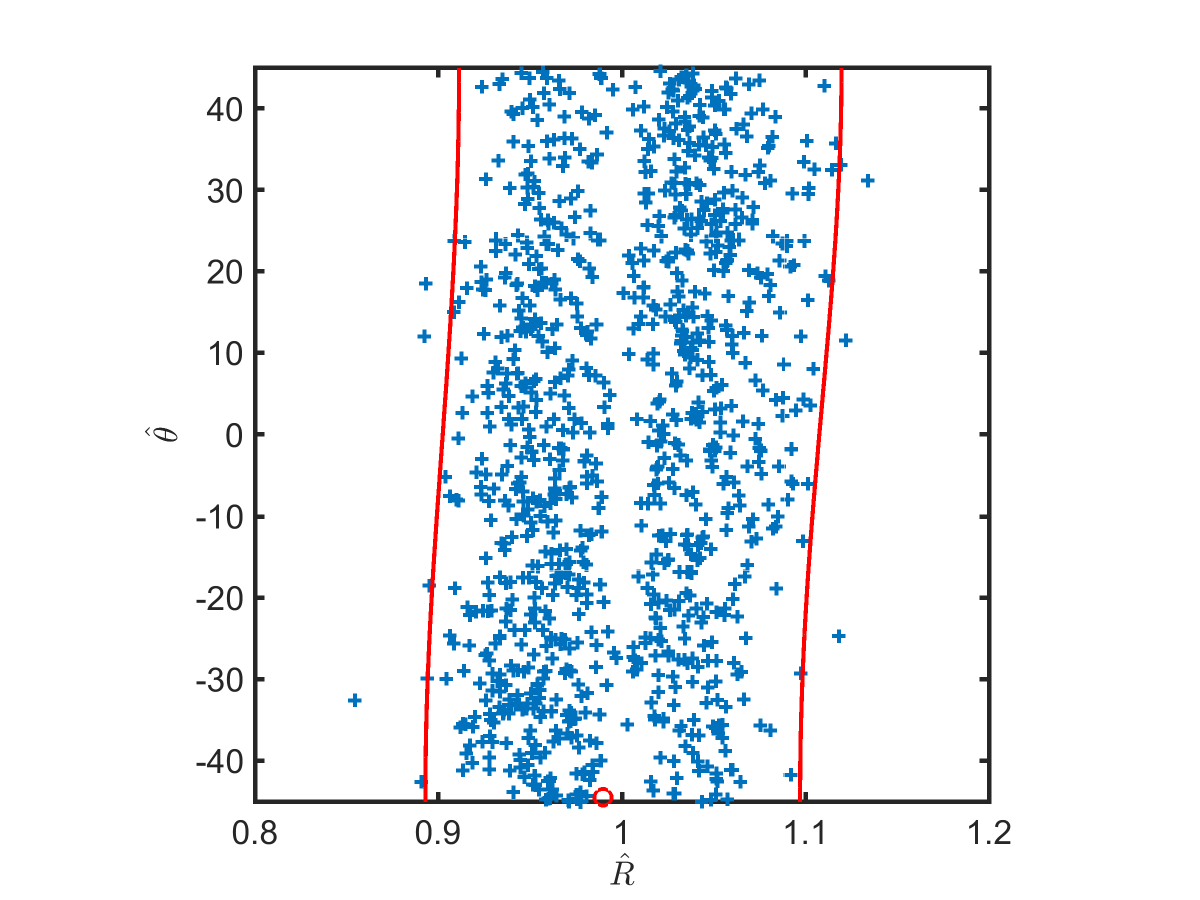}} %
	\subfloat[Mat\'ern, $L=800$, $N=1024$]{\label{subfig:M2} \includegraphics[width=0.33\linewidth, trim = 35pt 0pt 75pt 0pt, clip, clip]{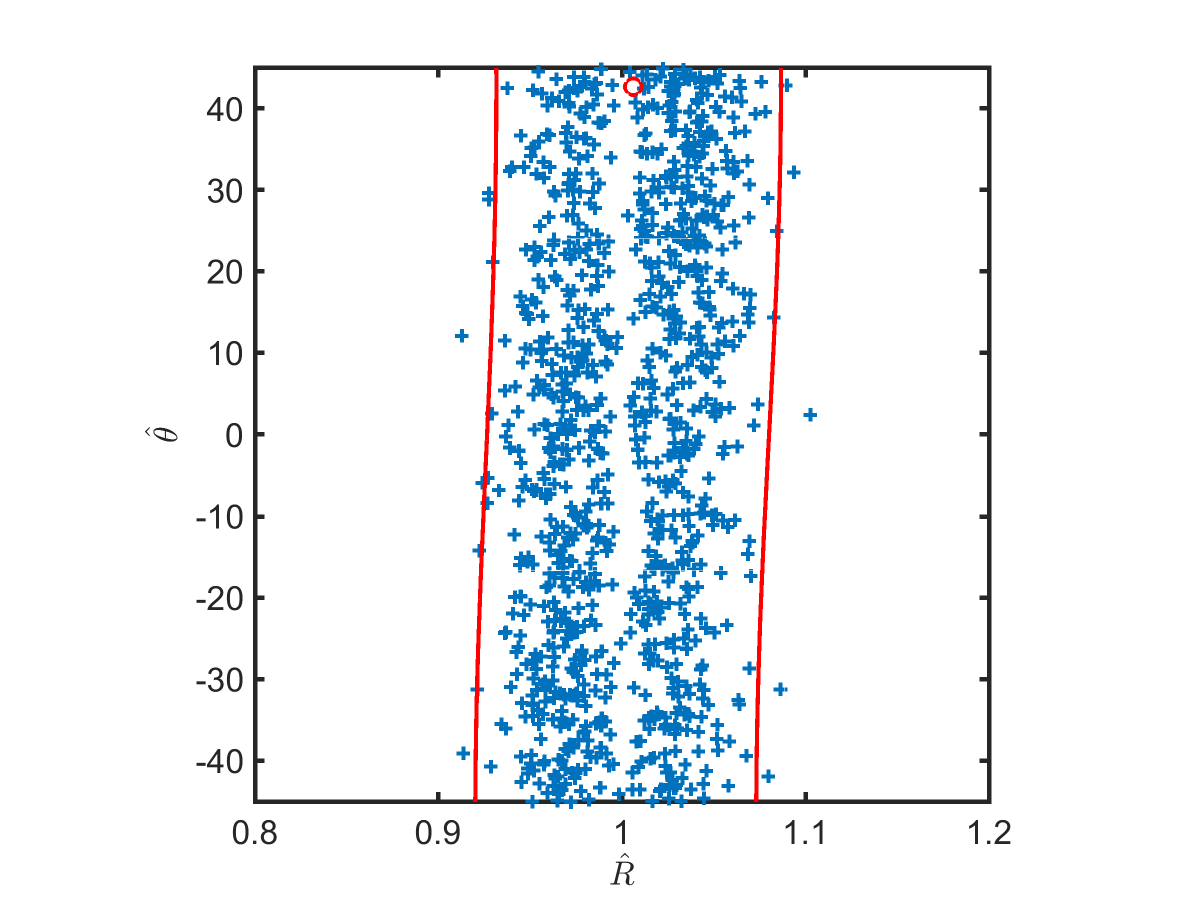}} %
	\subfloat[Mat\'ern, $L=1000$, $N=1600$]{\label{subfig:M3} \includegraphics[width=0.33\linewidth,trim = 35pt 0pt 75pt 0pt, clip, clip]{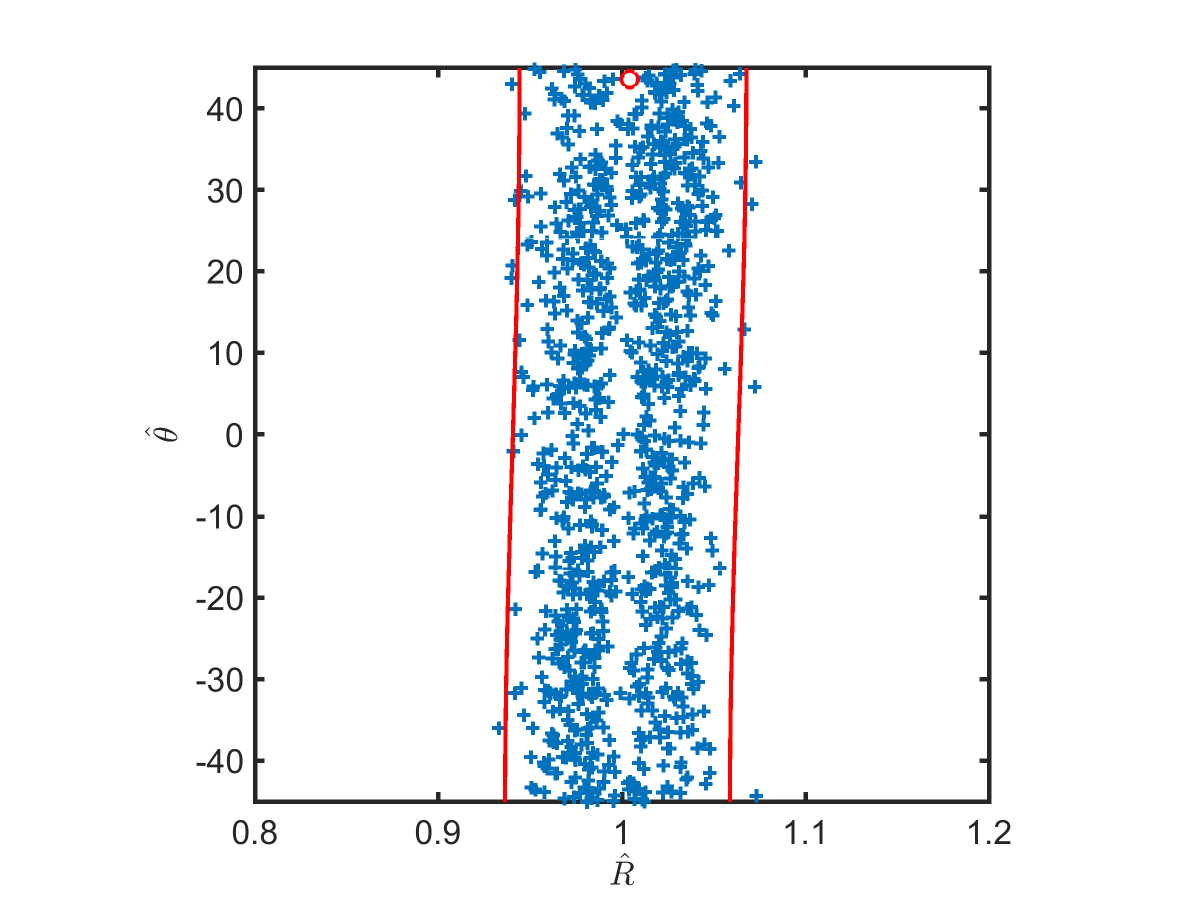}} %
	\caption{Non-parametric probability region estimation for isotropic scattered data. The initial lattice SRF is defined over a square lattice with side $L=600,800,1000$. The continuous curve corresponds to  $95\%$ non-parametric probability region calculated with anisotropy parameters estimated from $\overline{\bs{\cQ}}$. (a)--(c)~Anisotropy estimates (crosses) generated from $1000$ random samples obtained from a zero-mean, unit-variance isotropic Gaussian SRF with $\xi=28.3$. (d)--(f)~Anisotropy estimates (crosses) generated from $1000$ random scattered samples obtained from a zero-mean, unit-variance Mat\'ern SRF with $\xi=10$.}
	\label{fig:scattered_iso}
\end{figure*}

Figure~\ref{fig:scattered_iso} demonstrates the non-parametric probability regions and the scatter cloud of anisotropy estimates for scattered data sampled from isotropic Gaussian and Mat\'ern lattice SRFs of increasing side $L$. The absence of estimates near the $R = 1$ line agrees with the existence of a JPDF node at $R=1$ as discussed in Section~\ref{sec:isot-test}. For smaller domains the anisotropy estimates deviate from isotropy.

The computational complexity of natural neighbors interpolation is $\mathcal{O}((M+N)\log N)$ \citep{Park2006}, where $M$ is the number of the interpolation points. The complexity of derivative estimation using centered differences is $\mathcal{O}(M)$. Hence, $\mathcal{O}(2M)$ operations are needed for computing $\check{\partial_i}X(\bs{s})$, $i=1,2$ and $\mathcal{O}(3M)$ operations for $\hat{Q}_{ij}$. Thus the total computation time $t_{M,N}$ is of $\mathcal{O}(5M+(M+N)\log N)$, from which we obtain $t_{M',N'}/t_{M,N} = (5M'+(M'+N')\log N')/(5M+(M+N)\log N)$. For $M=200^2$, $t_{N=5184}/t_{N=576} = 1.28$. The time ratio obtained from Table~\ref{t:gaussian} is $\bar{t}_{N=5184}/\bar{t}_{N=576} = 1.36 \pm 0.06$ and from Table~\ref{t:matern} is $\bar{t}_{N=5184}/\bar{t}_{N=576} = 1.31 \pm 0.06$. For $M=100^2$, $t_{N=5184}/t_{N=576} = 1.53$ while the simulation times (average times per anisotropy estimation per processor core) obtained for 1000 realizations of isotropic Mat\'ern covariance (not shown here) are $\bar{t}_{N=576} = 36.4 \pm 1.5 \text{ msec}$ and $\bar{t}_{N=5184} = 58.3 \pm 2.6 \text{ msec}$, giving $\bar{t}_{N=5184}/\bar{t}_{N=576} = 1.60 \pm 0.10$.

\subsection{Case Study: Radiation Exposure}
\label{ssec:simulations-real}

We study anisotropy in two data sets of daily averages of radioactivity gamma dose rates over part of the Federal Republic of Germany, which was provided by the German automatic radioactivity monitoring network for the Spatial Interpolation Comparison (SIC 2004) exercise~\citep{dubois2005spatial}. Dose rates are measured in nanosieverts per hour (nSv/h). The \textit{background} data set corresponds to typical radioactivity measurements ($\approx 100$~nSv/h), which follow the Gaussian distribution (graph not shown), and thus their skewness and excess kurtosis coefficients are close to zero. The \textit{emergency} data includes a simulated local release of radioactivity  which results in five dose rate ``measurements'' around  $10$ times above background (exceeding $1000$~nSv/h). These measurements are aligned in the East-West direction. Table~\ref{t:SIC2004} summarizes the statistics of both data sets. The two rightmost columns show the CHI-based estimates of anisotropy parameters.
\begin{table*}
    \renewcommand{\arraystretch}{1.3}
    \setlength{\tabcolsep}{3pt}
    \centering
    \caption{Summary statistics of radioactivity dose rate exhaustive data sets (units are in nanosieverts per hour) and CHI anisotropy estimates. Abbreviations: min: minimum sample value; med: median sample value, max: maximum sample value; std: sample standard deviation; skew: sample skewness coefficient; kurt: sample excess kurtosis coefficient; $\cR, \,\ctheta$: estimates of anisotropy parameters.}
    \label{t:SIC2004}
    \begin{tabular}{l c c c c c c c c c}
    	\hline\hline
    	$N = 1008$ & min  & mean  & med  & max    & std  & skew & kurt  & $\cR$ & $\ctheta$     \\
    	\hline
    	Background & 57.0 & 97.7  & 98.6 & 180.0  & 19.6 & 0.4  & 0.6   & 1.18  & $7.36^\circ$  \\
    	Emergency  & 57.0 & 106.1 & 98.9 & 1528.2 & 92.5 & 11.3 & 144.1 & 0.45  & $-0.75^\circ$ \\
    	\hline\hline
    \end{tabular}
\end{table*}
Since the 95\% confidence interval for isotropy is $(R_{-}, R_{+}) = (0.92, 1.08)$, this dataset can be considered as slightly anisotropic. The direction of anisotropy is different in the two sets: In the background set the axis $A_{1}$ is tilted with respect to the x-axis (which is aligned with the E-W direction) by $7.36^\circ$, while the dominant anisotropy axis is $A_{2}$ since $\xi_{2} =1.18 \xi_{1}$. $A_{2}$ is closer to the y-axis, implying a dominant North-South anisotropy. In the emergency set the axis $A_{1}$ is slightly tilted with respect to the x-axis (by $-0.75^\circ$), and the dominant anisotropy axis is $A_{1}$ since $\xi_{2} = 0.45 \xi_{1}$. Since $A_{1}$ is closer to the x-axis, this implies that the radioactive plume reverses the
dominant  anisotropy direction to East-West.

\begin{figure}
    \centering
    \includegraphics[width=0.9\linewidth, trim = 0pt 0pt 0pt 100pt, clip]{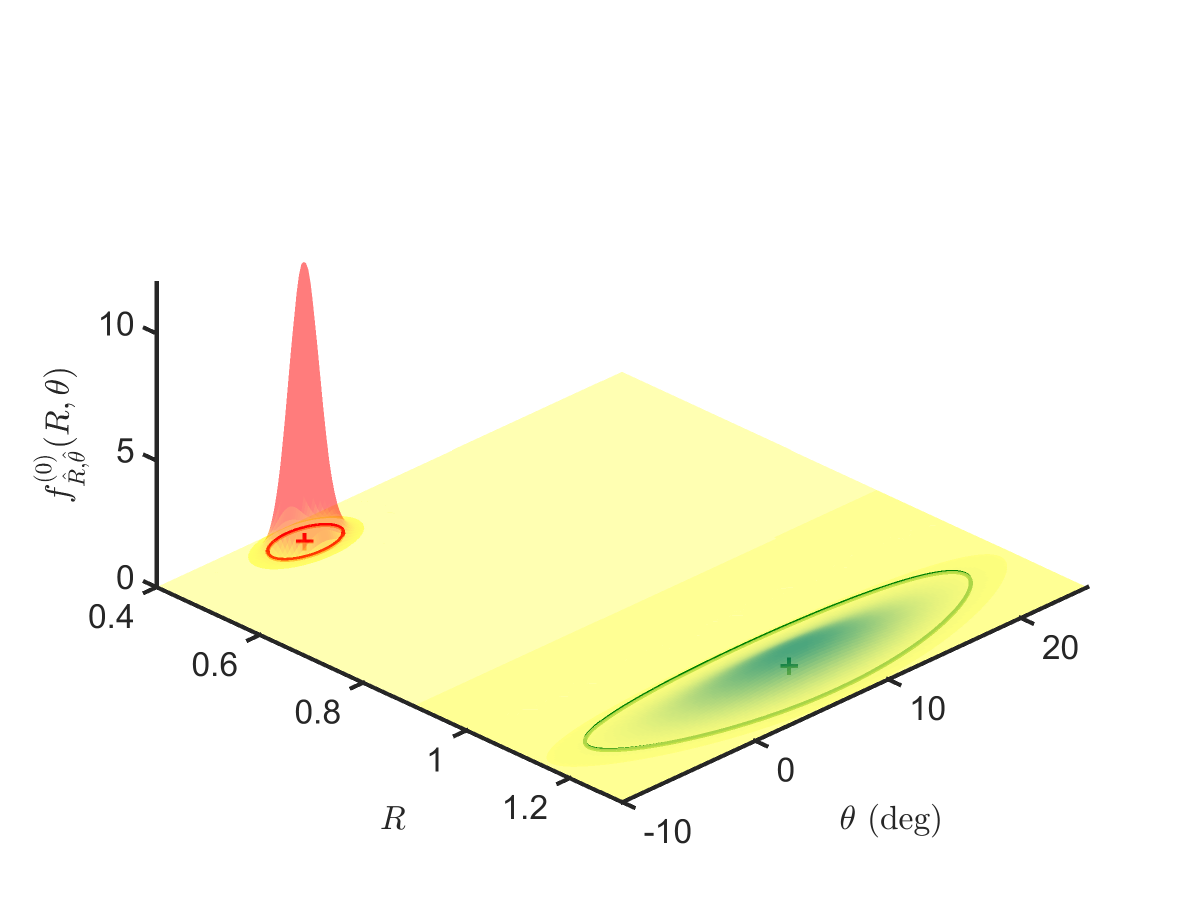}
    \caption{Non-parametric joint PDF shaded surfaces (red and green online) and 95\% confidence regions (solid contours) for the radioactivity dose rate data sets: background measurements (right) and emergency data simulation (left).}
    \label{fig:SIC2004}
\end{figure}

We calculate the non-parametric JPDF and the 95\% confidence regions of the anisotropy statistics based on the estimated anisotropy parameters (i.e., by CHI anisotropy estimation on gridded values obtained by natural neighbor interpolation) for both sets. The results are shown in Fig.~\ref{fig:SIC2004}.
There is no overlap of the two joint density functions, and the contours corresponding to the 95\%
confidence regions do not intersect. These patterns suggest statistically significant  anisotropy difference between the background and the emergency data due to the elevated values of the dose rate in the East-West direction of the spreading plume which changes the orientation of the major anisotropy axis.


\section{Discussion and Conclusions}
\label{sec:concl}
This work focuses on the estimation of geometric anisotropy in scattered or grid-based two-dimensional data.
We derive explicit expressions for the joint PDF of the anisotropy statistics, given by equations~\eqref{eq:triv-pdf-scaled-approx}-\eqref{eq:f-anispar}, and for the corresponding anisotropy probability regions at any level, i.e., equation~\eqref{eq:confi-reg}. The main assumptions used are that (i)~the data are drawn from a jointly Gaussian, stationary and differentiable random field and (ii)~the covariance function is short-ranged.

We also derive a non-parametric approximation for the joint PDF of the anisotropy statistics, which can be used if the covariance function is unknown \textit{a priori}, or if estimation of the covariance is not desired. The non-parametric approximation of the anisotropy joint PDF is given by~\eqref{eq:f-anis_np_both}. The corresponding equation for the non-parametric approximation of the probability region is given by~\eqref{eq:cr-R-theta}. We also derive probability intervals for the anisotropy ratio under the hypothesis that the sample comes from an isotropic random field. These probability intervals are used to formulate a non-parametric test of the isotropic hypothesis. We illustrate the application of the joint PDF and the probability regions with simulated and real data.

The Gaussian assumption is used in the decomposition of the covariance matrices $\CQQ$ and $\CQQ^{(0)}$, i.e., to derive equations~\eqref{eq:Cijkl} and~\eqref{eq:Cvv-nonpar} by means of the Wick-Isserlis theorem. This decomposition can be justified in principle, albeit approximately, even for non-Gaussian densities, based on
optimal variational (Gaussian) approximations. Practical application of the derived formulas requires the estimation of anisotropy statistics using CHI. Accurate estimation based on CHI requires in addition to (i) and (ii) above the following: (iii) a large sample size,  $N \to \infty$ and (iv) a sample domain that is large with respect to the correlation area. The latter is difficult to satisfy in case of large anisotropy $(\eR \gg 1$ or $\eR \ll 1)$. In such cases, the CHI estimate tends to underestimate the actual anisotropy. CHI anisotropy estimates include biases due to (v) the finite step size of the grid and (vi) interpolation (in the case of scattered data).

In summary, our approach consists of the following steps: (i) If necessary, preprocess the data in order to remove trends and use transforms to reduce deviations from the Gaussian distribution (e.g., Box-Cox transform).	(ii) Choose an interpolator which provides smooth interpolation surfaces. Several interpolators were tested and compared in~\citep{dth08}. The interpolation grid should be dense to allow accurate approximation of the spatial derivatives. (iii) Compute the partial derivatives of the interpolated surface. (iv) Compute $\hat{Q}_{ij}$ and use Theorem~\ref{theor:aniso} to obtain anisotropy parameter estimates. (v) The non-parametric JPDF is obtained from Theorem~\ref{theor:JPDF-NP}. The probability regions are provided by Lemma~\ref{theor:CR} using the approximate non-parametric covariance matrix, i.e.~\eqref{eq:Cvv-nonpar}. (vi) The isotropy test of Theorem~\ref{theor:isot-test-1} can be used to test for the presence of anisotropy. (vii) The differentiability assumption can be tested \emph{a posteriori} by determining the optimal anisotropic variogram model using standard geostatistical procedures. In this step, the CHI anisotropy estimates can be used to fix the anisotropy parameters or to provide informed initial guesses for likelihood optimization.

Our approach provides a computationally efficient, albeit approximate, method of geometric anisotropy estimation in two dimensions, because the analytical expressions derived above can be evaluated with minimal computational cost. The most computationally intensive part is the interpolation of scattered data onto regular grids in order to calculate derivatives. We use natural neighbor interpolation which is computationally fast (its complexity is essentially determined by Voronoi tesselation). For small datasets, the computation time scales linearly with the number of nodes $M$ of the interpolation grid, while for large datasets the computational cost is dominated by $\mathcal{O}(M \log N)$ where $N$ is the number of data points. A formulation of the natural neighbor interpolation algorithm which directly provides the partial derivatives of the interpolated surface is also available~\cite[Appendix A1]{Sambridge95}.

Our approach could be useful in estimating anisotropy in big data sets. In addition, the non-parametric JPDF can be used as an anisotropy prior in Bayesian and copula analyses~\citep{Kazianka13}. The method also provides initial estimates for maximum likelihood estimation of spatial anisotropic models~\citep{INTAMAP}. Furthermore, it furnishes an easily computable indicator of physical change in spatially extended systems based on the comparison of anisotropy probability regions.

Straightforward extension of this work is possible for the joint lognormal distribution along the lines of~\citep{dth08}. The global statistical measures of anisotropy can be efficiently calculated for large domains and can thus provide a useful statistic for large data sets. Local variations  of anisotropy can also be investigated using windowing methods. Capturing such local variability has applications in the analysis of medical images, e.g.~\citep{Richard10}. Currently, the solution of the non-linear CHI equations  for $d>2$ is not available in closed form. Hence, an analytical expression of the anisotropy joint PDF in higher than two dimensions is not yet feasible. Another path for future research is the development of an anisotropy detection method which will involve local integrals of the field values. Such an approach, if analytically tractable, will
apply to non-differentiable random fields as well.

\section*{Acknowledgment}
This work was funded by the European Commission, under the 6th FP, by the Contract N. 033811 with the DG INFSO, action line IST-2005-2.5.12 ICT for Environmental Risk Management. The views expressed herein are those of the authors and not necessarily of the European Commission.

We would like to thank Prof. Athanasios Liavas (School of Electronic and Computer Engineering, Technical University of Crete) for a careful reading of the manuscript and for suggesting improvements. In addition, we thank two anonymous reviewers for their valuable input.


\section*{Appendix A: Proof of Lemma~\ref{theor:CQQ}}
\label{app:A}
\renewcommand{\theequation}{A-\arabic{equation}}
\renewcommand{\thesubsection}{A-\arabic{subsection}}
\setcounter{equation}{0}

\begin{proof}
    Using the definition~\eqref{eq:Qdef} we obtain
    \begin{align}
        \label{eq:covm_def}
        C_{ij;kl}  &= \cov{\frac{1}{N}\sum_{n=1}^N {X}_{ij}(\bs{s}_n),
        \frac{1}{N}\sum_{m=1}^N {X}_{kl}(\bs{s}_m)} \nonumber \\
    &= \frac{1}{N^2} \sum_{n,m} \cov{{X}_{ij}(\bs{s}_n),{X}_{kl}(\bs{s}_m)}.
    \end{align}

    Due to the stationarity of $X(\bs{s})$, the double series in~\eqref{eq:covm_def} is reduced to
    a single series over all $(N^2)$ lag vectors $\bs{r}_{nm}=\bs{s}_n - \bs{s}_m$ $(n,m=1,\ldots,N)$, i.e.,
    \begin{align}
        \label{eq:cmexpansion}
        C_{ij;kl} & = \frac{1}{N^2} \sum_{\bs{r}_{nm}} \cov{X_{ij}(\bs{s}_0),X_{kl}(\bs{s}_0+\bs{r}_{nm})} \nonumber \\
                  & = \frac{1}{N} \cov{X_{ij}(\bs{0}),X_{kl}(\bs{0})} + \frac{1}{N^2} \sum_{\bs{r}_{nm} \neq \bs{0}} \cov{X_{ij}(\bs{0}),X_{kl}(\bs{r}_{nm})}.
    \end{align}
    \textit{Covariance of the gradient tensor:}
    Let $\bs{r}$ denote any lag vector (including $\bs{r}=\bs{0}$) between two points. Based on the definition of the covariance function it follows that
    \begin{equation}
        \label{eq:prodcov}
        \cov{X_{ij}(\bs{0}),X_{kl}(\bs{r})} = \E{X_{ij}(\bs{0}) X_{kl}(\bs{r})} -\E{X_{ij}(\bs{0})} \E{X_{kl}(\bs{r})}.
    \end{equation}

    Note that
    \begin{equation*}
    \E{X_{ij}(\bs{0}) X_{kl}(\bs{r})} = \E{\partial_{i}{X}(\bs{0}) \, \partial_{j}{X}(\bs{0}) \,
    \partial_{k}{X}(\bs{r}) \, \partial_{l}{X}(\bs{r})}.
    \end{equation*}

    For a  differentiable and stationary SRF $X(\bfs)$, the gradient component ${\partial}_{i}{X}(\bfs)$
    is a zero-mean Gaussian SRF with covariance function given by~\citep{Abrahamsen,Yaglom87}
    \begin{equation}
        \label{eq:cov-Xi}
        \E{\partial_{i}{X}(\bfs)\,\partial_{j}{X}(\bfs+{\bf r})}=
        - \frac{\partial^{2} \cxx(\bfr)}{\partial r_{i}\partial r_{j}}.
    \end{equation}

    Hence, $\E{X_{ij}(\bs{0}) \, X_{kl}(\bs{r})}$ can be calculated using the moment factorization property of multivariate normal distributions~\citep{Isserlis18,Wick50}
    \begin{align}
        \label{eq:EXijXkl}
        \E{X_{ij} (\bs{0}) X_{kl}(\bs{r})} &=
        \E{\partial_{i}{X}(\bs{0}) \partial_{j}{X}(\bs{0})} \E{\partial_{k}{X}(\bs{r}) \partial_{l}{X}(\bs{r})} \nonumber \\
        & + \E{\partial_{i}{X}(\bs{0}) \partial_{k}{X}(\bs{r})} \E{\partial_{j}{X}(\bs{0}) \partial_{l}{X}(\bs{r})}  \nonumber \\
        & + \E{\partial_{i}{X}(\bs{0}) \partial_{l}{X}(\bs{r})} \E{\partial_{j}{X}(\bs{0}) \partial_{k}{X}(\bs{r})}  \nonumber \\
        & = H_{ij}(\bs{0}) H_{kl}(\bs{0}) + H_{ik}(\bs{r}) H_{jl}(\bs{r}) + H_{il}(\bs{r}) H_{jk}(\bs{r}).
    \end{align}
    The last equality follows from Eq.~\eqref{eq:cov-Xi} and the definition~\eqref{eq:H} of CHM. The second term on the right-hand side of~\eqref{eq:prodcov} is   
    \begin{equation}
    \label{eq:EXijEXkl}
        \E{X_{ij}(\bs{0})}\E{X_{kl}(\bs{r})}  = \E{\partial_{i}{X}(\bs{0}) \partial_{j}{X}(\bs{0})}
        \E{\partial_{k}{X}(\bs{r}) \partial_{l}{X}(\bs{r})}
        = H_{ij}(\bs{0}) H_{kl}(\bs{0}).
    \end{equation}
    Thus, in light of~\eqref{eq:EXijXkl} and~\eqref{eq:EXijEXkl}, equation~\eqref{eq:prodcov} becomes
    \begin{equation}
    \label{eq:covXijXkl}
        \cov{X_{ij}(\bs{0}),X_{kl}(\bs{r})} = H_{ik}(\bs{r}) H_{jl}(\bs{r})
        + H_{il}(\bs{r}) H_{jk}(\bs{r}).
    \end{equation}
    Equation~\eqref{eq:Cijkl} follows from~\eqref{eq:cmexpansion}, \eqref{eq:covXijXkl}, and Theorem~\ref{theor:chi} for the zero-lag CHM.
\end{proof}


\section*{Appendix B: Proof of Lemma~\ref{theor:fQ}}
\label{app:B}
\renewcommand{\theequation}{B-\arabic{equation}}
\renewcommand{\thesubsection}{B-\arabic{subsection}}
\setcounter{equation}{0}

\begin{proof}
	To show that the JPDF of $\Qind$ tends asymptotically to the normal distribution,
	we use the multivariate CLT theorem. The classical CLT for scalar random variables is discussed in~\citep{GK:54,Levy:54,Feller:71}. The CLT  extension to vector random variables is as follows~\citep{Anderson03}:
	
	\emph{Assume $N$ independent and identically distributed vector variables $\bs{Z}_{k}$, $k=1,\dots,N$
	with mean $\bs{m}$ and covariance matrix $\bs{C}_{ZZ}$. Then, for $N \to \infty$ the joint distribution of the random vector $\bar{\bs{Z}}=(\bs{Z}_{1} + \dots + \bs{Z}_{N})/N$ tends to the multivariate normal distribution with mean $\bs{m}$ and covariance matrix $\bs{C}_{ZZ}/N$.}
	
	The above CLT is generalized to SRF averages. Loosely stated, an  average of a stationary  random field with finite-range correlations over $N \to \infty$ points tends to follow the joint normal probability distribution~\citep{Bouchaud}. Thus, the multivariate CLT applied to the random vector $\bs{Z}_{k} = \left( X_{11}(\bs{s}_k), X_{22}(\bs{s}_k), X_{12}(\bs{s}_k) \right)^{t}$
	leads to~\eqref{eq:f-v}.
	
	Next, we establish the condition for the SRFs to have finite correlation range. The $X_{ij}(\bs{s}_k)$ are stationary SRFs by virtue of the stationarity of  $X(\bfs)$. Hence, $\phi_{ijkl}(\bfr):= \cov{X_{ij}(\bs{s}),X_{kl}(\bs{s}+\bs{r})} = \cov{X_{ij}(\bs{0}),X_{kl}(\bs{r})}$. Using~\eqref{eq:covXijXkl}, $\phi_{ijkl}(\bs{r}) = H_{ik}(\bs{r})H_{jl}(\bs{r}) + H_{il}(\bs{r})H_{jk}(\bs{r})$. The correlation range of  $X_{ij}(\bs{s}_k)$ is determined by the integral	
	\begin{equation*}
	    V_{c} = \underset{i,j,k,l }{\max} \left( \frac{1}{\phi_{ijkl}(\bs{0})}\int {\mathrm d}\bs{r} \, \phi_{ijkl}(\bs{r}), \right).
	\end{equation*}
	
	Based on~\eqref{eq:covXijXkl}, $\phi_{ijkl}(\bs{0})= Q_{ij} \, Q_{kl} + Q_{il} \, Q_{jk}$ and thus $\phi_{ijkl}(\bs{0})$ has a finite value if $X(\bfs)$ has finite correlation lengths. We calculate  $\int_{\cal{D}} {\mathrm d}\bs{r} \, \phi_{ijkl}(\bs{r})$ in the asymptotic regime where $|\cal{D}|\to \infty$, and we express the integral in terms of the Fourier transform of $\cxx({\bf r})$. Any permissible covariance function $\cxx({\bf r})$, where ${\bf r} \in \mathbb{R}^2$, admits the following pair of transformations, where $\tilde{C}(\bs{k})$ is the \emph{spectral density}:	
	\begin{align*}
	    \cxx({\bf r}) = & \frac{1}{(2\pi)^2} \int {\mathrm d}\bs{k} \; \ee^{\jmath \bs{k} \cdot \bs{r}} \tilde{C}(\bs{k}),\\
	    \tilde{C}(\bs{k}) = & \int {\mathrm d}\bs{r} \, \ee^{- \jmath \bs{k} \cdot \bs{r}} \, \cxx({\bf r}).
	\end{align*}
	
	 \noindent Based on the above, it follows that $H_{ij}(\bs{r}) = (2\pi)^{-2} \int {\mathrm d}\bs{k} \, k_{i} k_{j} \, \ee^{\jmath \bs{k} \cdot \bs{r}} \tilde{C}(\bs{k})$, and thus	 
	\begin{equation*}
		\int {\mathrm d}\bs{r} \, \phi_{ijkl}(\bs{r}) = \frac{1}{(2\pi)^2} \int {\mathrm d}\bs{k} \, k_{i} k_{j} k_{k} k_{l} \, [ \tilde{C}(\bs{k}) ]^2.
	\end{equation*}
	
	\noindent In the above, $\jmath= \sqrt{-1}$, $ \bs{k} \cdot \bs{r} = k_{1}r_{1} + k_{2} r_{2}$ is the inner vector product, and $\int {\mathrm d}\bs{k}= \int_{-\infty}^{\infty}{\mathrm d}k_{1} \, \int_{-\infty}^{\infty}{\mathrm d}k_{2}$ or $\int {\mathrm d}\bs{k}= \int_{0}^{\infty} k\, {\mathrm d}k \, \int_{0}^{2\pi}{\mathrm d}\phi$ in polar coordinates. The existence of the above integral depends on the behavior of $\tilde{C}(\bs{k})$ at $\kk=0$ and $\kk \to \infty$. Since $\cxx({\bf r})$ is short-ranged,  $\int {\mathrm d}\bs{r}\, \cxx({\bf r})=  \tilde{C}(\bs{0})$ is finite, and thus the integrand is well-behaved at $\kk=0$. At $\kk \to \infty$, the integral converges (using polar coordinates) if $[ \tilde{C}(\bs{k}) ]^2 $ decays asymptotically faster than $\kk^{-6-2\epsilon}$, where $\epsilon>0$. This ensures that $\phi_{ijkl}(\bfr)$ is short-ranged.

\end{proof}


\section*{Appendix C: Proof of Lemma~\ref{theor:fq}}
\label{app:C}
\renewcommand{\theequation}{C-\arabic{equation}}
\renewcommand{\thesubsection}{C-\arabic{subsection}}
\setcounter{equation}{0}

\begin{proof}
	
	The probability transformation $\Q \mapsto \bq$ is performed as follows: Since $\dim(\bq) = 2 < \dim(\Q)=3$, we append to $\bq$ the dummy variable $u={Q}_{11} \ge 0$ and then integrate over  $u$. Using definitions~\eqref{eq:qdiag} and~\eqref{eq:qoff}, the absolute value of the Jacobian determinant for the transformation $({Q}_{11},{Q}_{22},{Q}_{12}) \mapsto (u,\qd, \qo)$ is
    \begin{equation*}
        \bs{J}_{\bs{q}} = \frac{\partial({Q}_{11}, {Q}_{22}, {Q}_{12})}
        {\partial(u,\qd,\qo)}
        \Rightarrow \abs{\det(\bs{J}_{\bs{q}})} = u^2.
    \end{equation*}

	\noindent The dummy variable $u$ is integrated, leading to
    \begin{equation}
    \label{eq:f-qd-qo-int}
        f_{\hbq}(\bq; \bs{\eQ},{\CQQ}) = \int_{0}^\infty {\mathrm d}u \,
        f_{\Qind}(u,\qo u,\qd u; \bs{\eQ},{\CQQ}) \, u^2.
    \end{equation}

	\noindent In terms of $\qd$  and $\qo$, the exponent of the PDF $f_{\Qind}(\cdot)$, given by~\eqref{eq:f-v}, becomes
    \begin{equation}
        \label{eq:quadr}
        (\Q - \bs{\eQ})^t \, {\CQQ^{-1}} \, (\Q - \bs{\eQ}) =  z_{1}^{2}(\bq ; \CQQ) \, u^2 +
        z_{2}(\bq ; \bs{\eQ}, \CQQ) \, u +  \la_{1}(\bs{\eQ}, {\CQQ}).
    \end{equation}

	\noindent By virtue of the above, ~\eqref{eq:f-qd-qo-int} is expressed as follows
    \begin{equation*}
        f_{\hbq}(\bq; \bs{\eQ},{\CQQ}) = \la_{2} \, \int_{0}^{\infty} {\mathrm d}u \,
        u^2 \, \ee^{-\frac{1}{2} \left[ (u\, z_{1})^2 +  u\, z_{2} + \la_{1} \right]}.
    \end{equation*}

	\noindent According to~\eqref{eq:A}, $ z_{1}^{2}>0$ because $\CQQ$ is a covariance matrix;
	hence $\CQQ$ as well as  $\CQQ^{-1}$ are positive definite. Thus, the Gaussian integral above exists and its value is given by~\eqref{eq:f-qd-qo}.
\end{proof}


\section*{Appendix D: Proof of Theorem~\ref{theor:jpdf}}
\label{app:D}
\renewcommand{\theequation}{D-\arabic{equation}}
\renewcommand{\thesubsection}{D-\arabic{subsection}}
\setcounter{equation}{0}

\begin{proof}
	
    Equation~\eqref{eq:f-anispar} follows from the transformation $(\qd,\qo) \mapsto (R,\theta)$ with Jacobian matrix $\bs{J}_{R, \theta}$. The transformed PDF is given by $f_{\hR,\htheta}(R,\theta;\bs{\eQ},{\CQQ}) = f_{\hbq}(\bq; \bs{\eQ},{\CQQ}) \,\abs{\det(\bs{J}_{R,\theta})}$,  where ${\det(\bs{J}_{R, \theta})}$ is given by
    \begin{equation}
        \label{eq:J}
        {\det(\bs{J}_{R, \theta})} =   \left|
                                        \begin{array}{cc}
                                         \frac{\partial \qd}{\partial R} &   \frac{\partial \qd}{\partial \theta} \\ [\medskipamount]
                                         \frac{\partial \qo}{\partial R}  & \frac{\partial \qo}{\partial \theta} \\
                                        \end{array}
                                      \right|
        =  \frac{2 R \, \left( R^2-1 \right) }{\left(R^2 \cos^2\theta + \sin^2 \theta \right)^3}.
    \end{equation}
	Restricting the parameter space to $R \in [0,\infty)$ and $\theta \in [-\pi/4,\pi/4)$, or equivalently $R \in [1,\infty)$ and $\theta \in [-\pi/2,\pi/2)$, the transformation $(\qd,\qo) \mapsto (R,\theta)$ is one-to-one except at $(1,0)$ in $(\qd, \qo)$-space, which is mapped onto the straight line  $R=1$ in the $(R,\theta)$-space, in which the Jacobian~\eqref{eq:J} vanishes. Finally, using Lemma~\ref{theor:fq}, $\JPDF$ is given by~\eqref{eq:f-anispar}.
\end{proof}


\bibliographystyle{model2-names}

\end{document}